\documentclass[12pt]{amsart}

\usepackage{amsthm}
\usepackage{amsmath}
\usepackage{amssymb}
\usepackage{bm}

\usepackage{hyperref}
\usepackage{url}
\usepackage{latexsym}
\usepackage{algorithm}
\usepackage{algorithmic}
\usepackage{tabulary}
\usepackage{multirow}
\usepackage{booktabs}
\usepackage{wrapfig}
\usepackage[dvipdfmx]{graphicx}

\newtheorem{thm}{Theorem}[section]
\newtheorem{prop}{Proposition}[section]
\newtheorem{lem}{Lemma}[section]

\newtheorem{prob}{Problem}[section]

\newtheorem{assum}{Assumption}[section]

\renewcommand{\theenumi}{\roman{enumi}}

\begin{document}
\makeatletter

\begin{center}
\large{\bf The Number of Steps Needed for Nonconvex Optimization of a Deep Learning Optimizer is a Rational Function of Batch Size}\\
\end{center}\vspace{3mm}

\begin{center}
\textsc{Hideaki Iiduka}\\
Department of Computer Science, 
Meiji University,
1-1-1 Higashimita, Tama-ku, Kawasaki-shi, Kanagawa 214-8571 Japan. 
(iiduka@cs.meiji.ac.jp)
\end{center}

\vspace{2mm}

\footnotesize{
\noindent\begin{minipage}{14cm}
{\bf Abstract:}
Recently, convergence as well as convergence rate analyses of deep learning optimizers for nonconvex optimization have been widely studied.
Meanwhile, numerical evaluations for the optimizers have precisely clarified the relationship between batch size and the number of steps needed for training deep neural networks. 
The main contribution of this paper is to show theoretically that the number of steps needed for nonconvex optimization of each of the optimizers can be expressed as a rational function of batch size. 
Having these rational functions leads to two particularly important facts, which were validated numerically in previous studies. 
The first fact is that there exists an optimal batch size such that the number of steps needed for nonconvex optimization is minimized.
This implies that using larger batch sizes than the optimal batch size does not decrease the number of steps needed for nonconvex optimization.
The second fact is that the optimal batch size depends on the optimizer.
In particular, it is shown theoretically that momentum and Adam-type optimizers can exploit larger optimal batches and further reduce the minimum number of steps needed for nonconvex optimization than can the stochastic gradient descent optimizer.
\end{minipage}
 \\[5mm]


\hbox to14cm{\hrulefill}\par


\section{Introduction}\label{sec:1}
One way to train deep neural networks is to find the model parameters of the deep neural networks that minimize loss functions called the expected risk and empirical risk using first-order optimization methods \cite[Section 4]{bottou}.
The simplest optimizer is stochastic gradient descent (SGD) \cite{robb1951,zinkevich2003,nem2009,gha2012,gha2013}.
There have been many deep learning optimizers to accelerate SGD, such as momentum methods \cite{polyak1964,nest1983}
and adaptive methods, e.g., Adaptive Gradient (AdaGrad) \cite{adagrad}, Root Mean Square Propagation (RMSProp) \cite{rmsprop},
Adaptive Moment Estimation (Adam) \cite{adam}, and Adaptive Mean Square Gradient (AMSGrad) \cite{reddi2018} (Table 2 in \cite{Schmidt2021} lists useful deep learning optimizers).

Convergence and convergence rate analyses of deep learning optimizers have been widely studied for convex optimization \cite{zin2010,adam,reddi2018,luo2019,dun2020}.
Meanwhile, theoretical investigation of deep learning optimizers for nonconvex optimization is needed so that these optimizers can put into practice for nonconvex optimization in deep learning \cite{kxu2015,ar2017,vas2017}.

Convergence analyses of SGD for nonconvex optimization were studied in 
\cite{feh2020,chen2020,sca2020,loizou2021} (see \cite{gower2021,loizou2021} for convergence analyses of SGD for two classes of nonconvex optimization problems, quasar-convex and Polyak--Lojasiewicz optimization problems).
For example, Theorem 11 in \cite{sca2020} indicates that SGD with a diminishing learning rate $\alpha_k = 1/\sqrt{k}$ has $\mathcal{O}(1/\sqrt{K})$ convergence, where $K$ denotes the number of steps. 
Convergence analyses of SGD depending on the batch size were presented in \cite{chen2020}. 
In particular, Theorem 3.2 in \cite{chen2020} indicates that
running SGD with a diminishing learning rate $\alpha_k = 1/k$ and large batch size for sufficiently many steps
leads to convergence to a local minimizer of a sum of loss functions.
 
Convergence analyses of adaptive methods for nonconvex optimization were studied in \cite{spider,chen2019,ada,iiduka2021}. 
In \cite{chen2019}, it was shown that generalized Adam, which includes the Heavy-ball method, AdaGrad, RMSProp, AMSGrad, and AdaGrad with First Order Momentum (AdaFom), using a diminishing learning rate $\alpha_k = 1/\sqrt{k}$ has an $\mathcal{O}(\log K/\sqrt{K})$ convergence rate. 
AdaBelief (named for adapting stepsizes by the belief in observed gradients) using $\alpha_k = 1/\sqrt{k}$ has $\mathcal{O}(\log K/\sqrt{K})$ convergence \cite{ada}.
In \cite{iiduka2021}, a method was presented to unify useful adaptive methods such as AMSGrad and AdaBelief, and it was shown that the method with $\alpha_k = 1/\sqrt{k}$ has an $\mathcal{O}(1/\sqrt{K})$ convergence rate, which improves on the results in \cite{chen2019,ada}. 
A theoretical investigation of Stochastic Path-Integrated Differential EstimatoR (SPIDER) for $\epsilon$-approximation in nonconvex optimization was reported in \cite{spider}.
In particular, Theorem 2 in \cite{spider} clarified that SPIDER, which has a constant learning rate, for $\epsilon$-approximation must use the full-batch gradient with the number of samples $n$ or the stochastic gradient with batch size $\sqrt{n}$. 
 
Meanwhile, in \cite{shallue2019}, it was studied how increasing the batch size affects the performances of SGD, SGD with momentum \cite{polyak1964,momentum}, and Nesterov momentum \cite{nest1983,sut2013}. 
The relationships between batch size and performance for Adam and K-FAC (Kronecker-Factored Approximate Curvature \cite{kfac}) were studied in \cite{zhang2019}. 
In both studies, it was numerically shown that increasing batch size tends to decrease the number of steps $K$ needed for training deep neural networks,  
but with diminishing returns \cite[Figure 4]{shallue2019}, \cite[Figure 8]{zhang2019}. 
Moreover, it was shown that SGD with momentum and Nesterov momentum can exploit larger batches than SGD \cite[Figure 4]{shallue2019}, 
and that K-FAC and Adam can exploit larger batches than SGD with momentum \cite[Figure 5]{zhang2019}.
Thus, it was shown that momentum and adaptive methods can significantly reduce the number of steps $K$ needed for training deep neural networks \cite[Figure 4]{shallue2019}, \cite[Figure 5]{zhang2019}.

\subsection{Contribution}
The contribution of this paper is to construct a theory guaranteeing the useful numerical results in \cite{shallue2019,zhang2019}.
Table \ref{constant} (resp. Table \ref{diminishing}) summarizes our results for SGD, Nesterov momentum (N-Momentum), and Adam-type optimizers with a constant learning rate rule (resp. diminishing learning rate rule), described in Theorem \ref{thm:1} (resp. Theorem \ref{thm:2}). 
See Theorem \ref{thm:3} in Appendix for other result for the optimizers with a diminishing learning rate rule.
Figure \ref{constant_fig} (resp. Figure \ref{diminishing_fig}) visualizes the relationships between the optimizers for the results shown in Table \ref{constant} (resp. Table \ref{diminishing}) for an appropriately set momentum coefficient $\beta$.

\begin{table}[htpb]
\caption{Relationship between batch size $s$ and the number of steps $K_\epsilon$
needed for nonconvex optimization in the sense of (\ref{evaluation}) of optimizers with constant learning rates}\label{constant}
\begin{tabular}{l||c|c|c}
\bottomrule
\multirow{3}{*}{} & \multicolumn{3}{c}{Constant Learning Rate Rule} \\
                  & \multicolumn{3}{c}{($\alpha_k = \alpha \in (0,1]$, $\beta_k = \beta \in [0,b] \subset [0,1)$)}  \\ \cline{2-4}
                  & Rational Function  & Optimal Batch Size $s^\star$ & Minimum Steps $K_\epsilon(s^\star)$ \\
\hline
SGD
& $\displaystyle{K_\epsilon = \frac{A_\alpha s^2}{\epsilon^2 s - B_\alpha}}$ 
& $\displaystyle{\frac{d D L^2 n^2 \alpha}{\epsilon^2}}$      
& $\displaystyle{\frac{(dDLn)^2}{\epsilon^4}}$     \\ 
N-Momentum  
& $\displaystyle{K_\epsilon = \frac{A_\alpha s^2}{(\epsilon^2 - C_\beta) s - B_\alpha}}$      
& $\displaystyle{\frac{d D L^2 n^2 \alpha}{\tilde{b}\epsilon^2 - dDLn\beta}}$      
& $\displaystyle{\frac{(dDLn)^2 }{(\tilde{b}\epsilon^2 - dDLn \beta)^2}}$\\
Adam-type         
& $\displaystyle{K_\epsilon = \frac{A_\alpha s^2}{(\epsilon^2 - C_\beta) s - B_\alpha}}$      
& $\displaystyle{\frac{d D L^2 n^2 \alpha}{\tilde{\gamma}^2(\tilde{b}\epsilon^2 - dDLn\beta)h_0^*}}$      
& $\displaystyle{\frac{(dDLn)^2 H}{\tilde{\gamma}^2(\tilde{b}\epsilon^2 - dDLn\beta)^2 h_0^*}}$   \\
\bottomrule 
\end{tabular}
Let $\epsilon > 0$, $\tilde{b} := 1 -b$, $\tilde{\gamma} := 1 -\gamma$ ($\gamma \in [0,1)$),
and $H \geq h_0^* > 0$.
The number of samples is denoted by $n$,  
$\nabla f_i \colon \mathbb{R}^d \to \mathbb{R}^d$ ($i\in [n] := \{1,2,\ldots,n\}$) is Lipschitz continuous with Lipschitz constant $L_i$, and $L$ denotes the maximum value of $L_i$. 
$D$ is the upper bound of $(x_{k,i} - x_i)^2$ ($(x_i) \in \mathbb{R}^d$), where $(\bm{x}_k)_{k\in\mathbb{N}} = ((x_{k,i}))_{k\in\mathbb{N}}$ is generated by an optimizer. 
$A_\alpha$ and $B_\alpha$ are positive constants depending on a learning rate $\alpha$
and $C_\beta$ is a positive constant depending on a momentum coefficient $\beta$ (see Theorem \ref{thm:1} for detailed definitions of the constants).
\end{table}

\begin{table}[htpb]
\caption{Relationship between batch size $s$ and the number of steps $K_\epsilon$
needed for nonconvex optimization in the sense of (\ref{evaluation}) of optimizers with diminishing learning rates}\label{diminishing}
\begin{tabular}{l||c|c|c}
\bottomrule
\multirow{3}{*}{} & \multicolumn{3}{c}{Diminishing Learning Rate Rule} \\
                  & \multicolumn{3}{c}{($\alpha_k = \frac{\alpha}{\sqrt{k}}$, $\beta_k = \beta \in [0,b] \subset [0,1)$)}  \\ \cline{2-4}
                  & Rational Function  & Optimal Batch Size $s^\star$ & Minimum Steps $K_\epsilon(s^\star)$ \\
\hline
SGD
& $\displaystyle{K_\epsilon = \left\{\frac{A_\alpha s^2 + B_\alpha}{\epsilon^2 s}\right\}^2}$ 
& $\displaystyle{\sqrt{2} Ln \alpha}$      
& $\displaystyle{\frac{2(dDLn)^2}{\epsilon^4}}$     \\ 
N-Momentum  
& $\displaystyle{K_\epsilon = \left\{\frac{A_\alpha s^2 + B_\alpha}{(\epsilon^2 - C_\beta) s}\right\}^2}$      
& $\displaystyle{\sqrt{2} Ln \alpha}$      
& $\displaystyle{\frac{2(dDLn)^2}{(\tilde{b}\epsilon^2 - dDLn\beta)^2}}$\\
Adam-type         
& $\displaystyle{K_\epsilon = \left\{\frac{A_\alpha s^2 + B_\alpha}{(\epsilon^2 - C_\beta) s}\right\}^2}$      
& $\displaystyle{\frac{\sqrt{2} Ln \alpha}{\tilde{\gamma} \sqrt{Hh_0^*}}}$      
& $\displaystyle{\frac{2(dDLn)^2 H}{\tilde{\gamma}^2(\tilde{b}\epsilon^2 - dDLn \beta)^2 h_0^*}}$   \\
\bottomrule 
\end{tabular}
See Table \ref{constant} and Theorem \ref{thm:2} for definitions of the constants.
\end{table}

The main contribution of this paper is to clarify that  
\begin{itemize}
\item 
the number of steps $K = K_\epsilon$ needed for nonconvex optimization in the sense of 
\footnote{Jensen's inequality guarantees that (\ref{evaluation}) implies that 
$\min_{k\in [K]} \mathbb{E}\left[ \| \nabla f(\bm{x}_k)  \| \right] \leq \epsilon$.} 
\begin{align}\label{evaluation}
\min_{k\in [K]} \mathbb{E}\left[ \| \nabla f(\bm{x}_k)\|^2 \right] \leq \epsilon^2,
\end{align}
where $\nabla f \colon \mathbb{R}^d \to \mathbb{R}^d$ denotes the gradient of a nonconvex loss function $f \colon \mathbb{R}^d \to \mathbb{R}$, $\epsilon > 0$ is a precision accuracy, and the sequence $(\bm{x}_k)_{k\in\mathbb{N}} \subset \mathbb{R}^d$ is generated by a particular optimizer, such as one of SGD, N-Momentum, and Adam-type optimizers, can be expressed as a rational function of batch size $s$ 
(see the ``Rational Function" columns of Tables \ref{constant} and \ref{diminishing}).
\end{itemize}

\begin{figure}[htbp]
\begin{tabular}{cc}
\begin{minipage}[t]{0.5\hsize}
\centering
\includegraphics[width=1\textwidth]{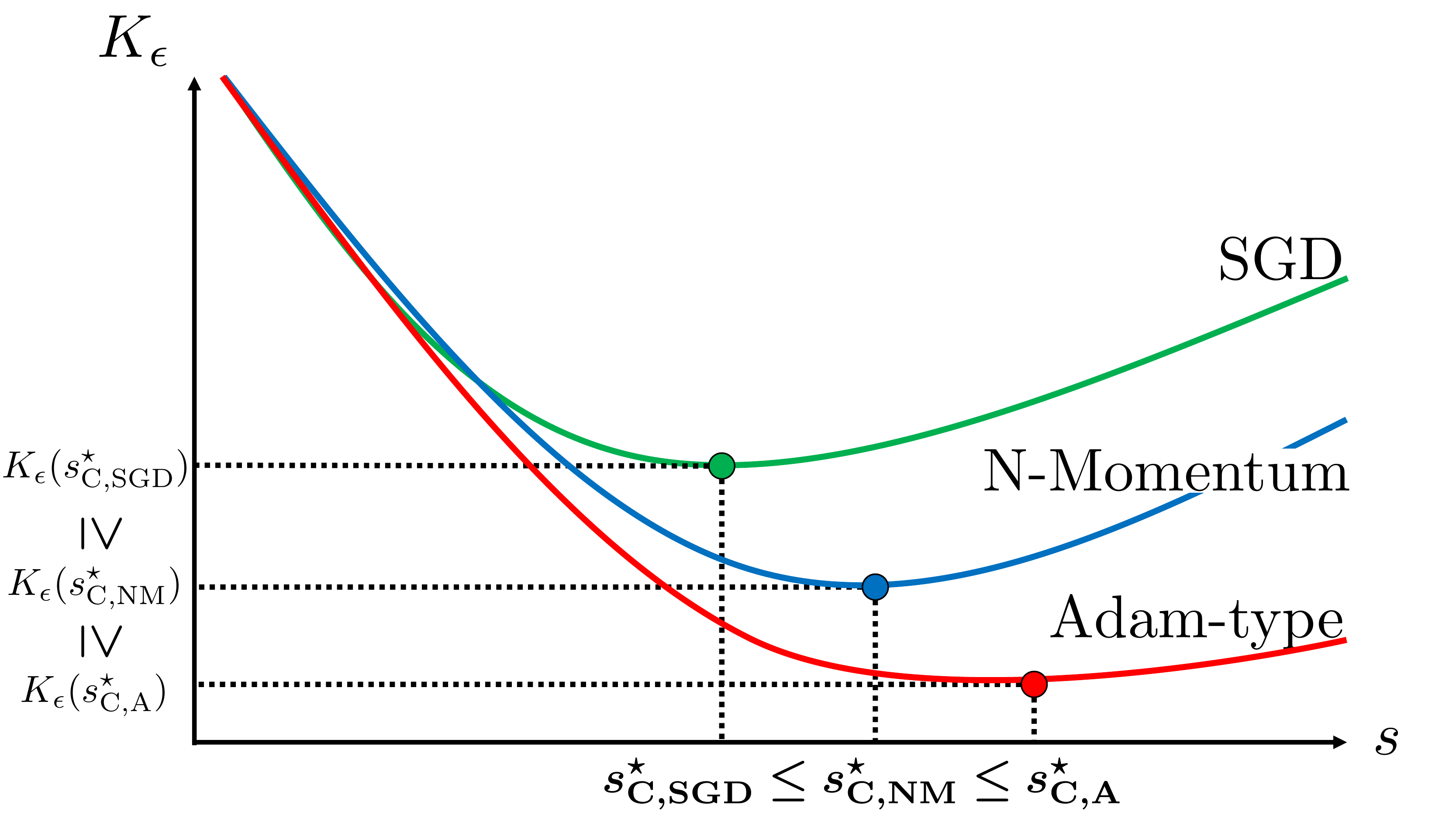}
\caption{Relationships between the optimizers in terms of the results in Table \ref{constant} 
(relations $s_{\mathrm{C,SGD}}^\star \leq s_{\mathrm{C,NM}}^\star \leq s_{\mathrm{C,A}}^\star$  hold generally, but those of  $K_\epsilon(s_{\mathrm{C,SGD}}^\star) \geq K_\epsilon(s_{\mathrm{C,NM}}^\star) \geq K_\epsilon(s_{\mathrm{C,A}}^\star)$ depend on momentum coefficient $\beta$)}
\label{constant_fig}
\end{minipage} &
\begin{minipage}[t]{0.5\hsize}
\centering
\includegraphics[width=1\textwidth]{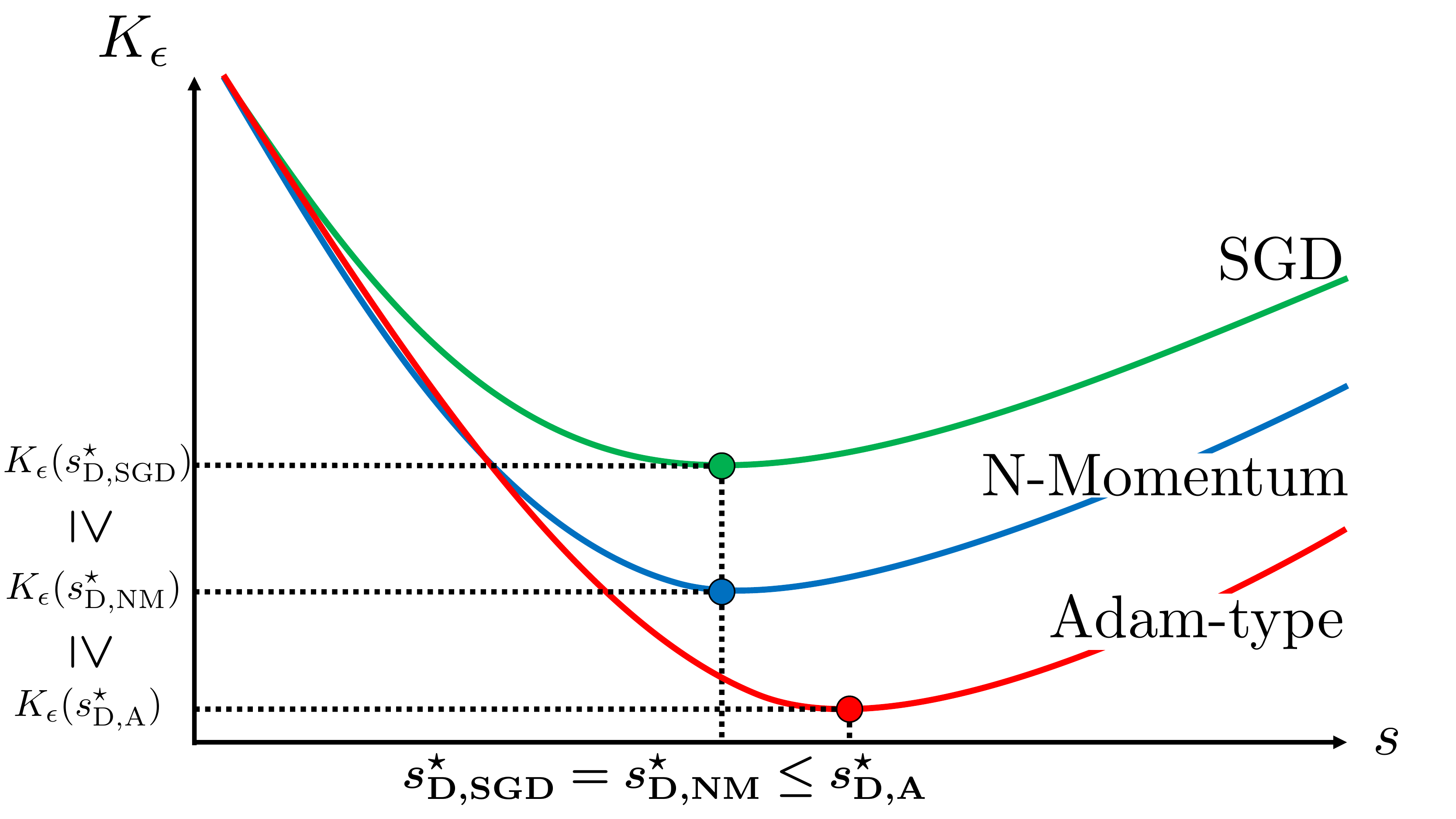}
\caption{Relationships between the optimizers in terms of the results in Table \ref{diminishing}
(relations $s_{\mathrm{D,SGD}}^\star = s_{\mathrm{D,NM}}^\star \leq s_{\mathrm{D,A}}^\star$  hold generally, but those of $K_\epsilon(s_{\mathrm{D,SGD}}^\star) \geq K_\epsilon(s_{\mathrm{D,NM}}^\star) \geq K_\epsilon(s_{\mathrm{D,A}}^\star)$ depend on momentum coefficient $\beta$)}
\label{diminishing_fig}
\end{minipage}
 \end{tabular}
\end{figure}      

The explicit forms of the rational functions imply the following two significant facts:
\begin{enumerate}
\item[(I)]
There exists an optimal batch size $s^\star$ such that $K_\epsilon(s)$ is minimized; specifically, $K_\epsilon(s)$ is monotone decreasing for $s \leq s^\star$ and  monotone increasing for $s \geq s^\star$. 
This fact guarantees theoretically the existences of the diminishing returns shown in \cite[Figure 4]{shallue2019}, \cite[Figure 8]{zhang2019}, which are such that increasing the batch size does not decrease the number of steps $K_\epsilon$.
\item[(II)]
The optimal batch size $s^\star$ and the minimum number of steps $K_\epsilon(s^\star)$ depend on the optimizer. In particular,
N-Momentum and Adam-type optimizers can exploit the same sized or larger batches ($s^\star$ in Tables \ref{constant} and \ref{diminishing} and Figures \ref{constant_fig} and \ref{diminishing_fig}) than can SGD. 
Furthermore, the dependence of N-Momentum and Adam-type optimizers on $\beta$ allows them to reduce the minimum number of steps ($K_\epsilon(s^\star)$ in Tables \ref{constant} and \ref{diminishing} and Figures \ref{constant_fig} and \ref{diminishing_fig}) more than can SGD (see Section \ref{sec:main} for details).
\end{enumerate}

\subsection*{Comparisons of Optimal Batch Sizes for Different Learning Rate Rules}
Tables \ref{constant} and \ref{diminishing} ensure that $K_\epsilon(s^\star)$ for Algorithm \ref{algo:1} using constant learning rates is almost the same as $K_\epsilon(s^\star)$ for Algorithm \ref{algo:1} using diminishing learning rates. 
Meanwhile, we would like to emphasize that the optimal batch size $s_{\mathrm{C}}^\star$ for Algorithm \ref{algo:1} using constant learning rates depends on $\epsilon$ and $\beta$, and the optimal batch size $s_{\mathrm{D}}^\star$ for Algorithm \ref{algo:1} using diminishing learning rates
does not depend on $\epsilon$ and $\beta$.
For example, under the precision accuracy $\epsilon = 10^{-1}$,
we can know the optimal batch sizes for N-Momentum with the frequently used parameter value $\alpha  = 10^{-3}$ are respectively 
\begin{align*}
s_{\mathrm{C,NM}}^\star = \frac{d D L^2 n^2 \epsilon^3}{\tilde{b}\epsilon^2 - dDLn\beta}
\text{ and }
s_{\mathrm{D,NM}}^\star = \sqrt{2}Ln \epsilon^3
\end{align*} 
before implementing N-Momentum.

\subsection{Notation}
$\mathbb{N}$ denotes the set of nonnegative integers. Let $n \in \mathbb{N} \backslash \{0\}$. We define $[n]:=\{1,2,\ldots,n\}$. 
$\mathbb{R}^d$ denotes $d$-dimensional Euclidean space with inner product $\langle \cdot, \cdot \rangle$ inducing norm $\| \cdot \|$. 
Let $\mathbb{S}_{++}^d$ be the set of $d \times d$ symmetric positive-definite matrices
and let $\mathbb{D}^d$ be the set of $d \times d$ diagonal matrices: $\mathbb{D}^d = \{ M \in \mathbb{R}^{d \times d} \colon M = \mathsf{diag}(x_i), \text{ } x_i \in \mathbb{R} \text{ } (i\in [d]) \}$.  
For a random variable $Z$, we use $\mathbb{E}[Z]$ to indicate its expectation.

\section{Nonconvex Optimization and Deep Learning Optimizers}
\subsection{Assumptions Regarding Loss Function and Gradient Estimation}\label{subsec:assumptions}
This paper considers optimization problems under the following assumptions. 

\begin{assum}\label{assum:0}
\text{ } 
{\em
\begin{enumerate} 
\item[(A1)][Loss function] $f_i \colon \mathbb{R}^d \to \mathbb{R}$ ($i \in [n]$) is differentiable and $f \colon \mathbb{R}^d \to \mathbb{R}$ is defined for all $\bm{x}\in \mathbb{R}^d$ by 
\begin{align*}
f(\bm{x}) := 
\frac{1}{n} \sum_{i=1}^n f_i (\bm{x}),
\end{align*}
where $n$ denotes the number of samples.
\item[(A2)][Gradient estimation] For each iteration $k$, optimizers sample a batch $\mathcal{S}_k \subset [n]$ of size $s := |\mathcal{S}_k|$ independently of $k$ and estimate the full gradient $\nabla f$ as 
\begin{align*}
\nabla f_{\mathcal{S}_k} := \frac{1}{s} \sum_{i\in \mathcal{S}_k} \nabla f_i.
\end{align*} 
\item[(A3)][Gradient boundedness] There exists a positive number $G$ such that, for all $\bm{x}\in X$, 
\begin{align}
\mathbb{E}\left[\|\nabla f_{\mathcal{S}_k} (\bm{x})\|^2 \right] \leq \frac{G^2}{s^2},
\end{align}
where $X$ is a subset of $\mathbb{R}^d$.
\end{enumerate}
}
\end{assum}

Assumption (A1) is a standard one for nonconvex optimization in deep neural networks (see, e.g., \cite[(2)]{chen2019} and \cite[(1.2)]{spider}). Assumption (A2) is needed for the optimizers to work (see, e.g., \cite[Section 2]{chen2019} and \cite[Notation section]{spider}).
Assumption (A3) is used to analyze the optimizers.
Assumption (A3) holds if each of the following holds (see Proposition \ref{prop:1} in Appendix for details):
\begin{enumerate}
\item[(G1)]
$X \subset \mathbb{R}^d$ is bounded, the gradient $\nabla f_i$ is Lipschitz continuous with Lipschitz constant $L_i$,
and $S_i := \{ \bm{x}^* \in \mathbb{R}^d \colon \nabla f_i(\bm{x}^*) = \bm{0} \} \neq \emptyset$ ($i\in [n]$), 
where 
$L := \max_{i\in [n]} L_i$.
(If we define $G_{k,L} := \sup_{\bm{x} \in X}  \sum_{i\in \mathcal{S}_k} \|\nabla f_i(\bm{x}) \|$, 
then we can take $G := \sup_{k\in \mathbb{N}} G_{k,L}$.)\item[(G2)]
$X \subset \mathbb{R}^d$ is bounded and closed. (If we define $G_{k} := \sup_{\bm{x} \in X}  \sum_{i\in \mathcal{S}_k} \|\nabla f_i(\bm{x}) \|$, then we can take $G := \sup_{k\in \mathbb{N}} G_{k}$.)\end{enumerate}

\subsection{Nonconvex Optimization in Deep Learning}
\label{sec:2}
This paper considers the following problem \cite{chen2019,ada}.

\begin{prob}\label{prob:1}
Under Assumption \ref{assum:0}, we would like to find a local minimizer $\bm{x}^\star$ of $f$ over $\mathbb{R}^d$, i.e.,
\begin{align*}
\bm{x}^\star \in X^\star := 
\left\{ \bm{x} \in \mathbb{R}^d \colon \nabla f(\bm{x}) = \bm{0} \right\}.
\end{align*}
\end{prob}

If $f$ is convex \cite{adam,reddi2018}, then the solution to Problem \ref{prob:1} is a global minimizer of $f$ over $\mathbb{R}^d$.
See the third and fourth paragraphs of Section \ref{sec:1} for the previous studies on Problem \ref{prob:1}.

\subsection{Deep Learning Optimizers}\label{subsec:optimizers}
There are many deep learning optimizers \cite[Table 2]{Schmidt2021}.
In this paper, we consider the following algorithm (Algorithm \ref{algo:1}),
which is a unified algorithm for useful optimizers, for example,  
N-Momentum \cite{nest1983,sut2013}, AMSGrad \cite{reddi2018,chen2019}, AMSBound \cite{luo2019}, and AdaBelief \cite{ada}, listed in Table \ref{table:ex} in Appendix.

\begin{algorithm} 
\caption{Deep learning optimizer for solving Problem \ref{prob:1}} 
\label{algo:1} 
\begin{algorithmic}[1] 
\REQUIRE
$(\alpha_k)_{k\in\mathbb{N}} \subset (0,1]$, $(\beta_k)_{k\in\mathbb{N}} \subset [0,b] \subset [0,1)$, 
${\gamma} \in [0,1)$
\STATE
$k \gets 0$, $\bm{x}_0, \bm{m}_{-1} := \bm{0} \in \mathbb{R}^d$, 
$\mathsf{H}_0 \in \mathbb{S}_{++}^d \cap \mathbb{D}^d$,
$\mathcal{S}_0 \subset [n]$
\LOOP 
\STATE 
$\bm{m}_k := \beta_k \bm{m}_{k-1} + (1-\beta_k) \nabla f_{\mathcal{S}_k}(\bm{x}_k)$
\STATE
$\displaystyle{\hat{\bm{m}}_k := \frac{\bm{m}_n}{1-{\gamma}^{k+1}}}$
\STATE
$\mathsf{H}_k \in \mathbb{S}_{++}^d \cap \mathbb{D}^d$ \text{ } (see Table \ref{table:ex} for examples of $\mathsf{H}_k$)
\STATE 
Find $\bm{\mathsf{d}}_k \in \mathbb{R}^d$ that solves $\mathsf{H}_k \bm{\mathsf{d}} = - \hat{\bm{m}}_k$
\STATE 
$\bm{x}_{k+1} := \bm{x}_k + \alpha_k \bm{\mathsf{d}}_k$
\STATE $k \gets k+1$
\ENDLOOP 
\end{algorithmic}
\end{algorithm}

The useful optimizers, such as N-Momentum, AMSGrad, AMSBound, and AdaBelief (Table \ref{table:ex}), all satisfy the following conditions:

\begin{assum}\label{assum:1}
{\em 
The sequence $(\mathsf{H}_k)_{k\in\mathbb{N}} \subset \mathbb{S}_{++}^d \cap \mathbb{D}^d$, with $\mathsf{H}_k := \mathsf{diag}(h_{k,i})$, in Algorithm \ref{algo:1} satisfies the following conditions:
\begin{enumerate}
\item[(A4)] $h_{k+1,i} \geq h_{k,i}$ almost surely for all $k\in\mathbb{N}$ and all $i\in [d]$;
\item[(A5)] For all $i\in [d]$, a positive number $H_i$ exists such that $\sup_{k \in \mathbb{N}} \mathbb{E}[h_{k,i}]  \leq H_i$.
\end{enumerate}
Moreover, the following condition holds:
\begin{enumerate}
\item[(A6)] $D := \max_{i\in [d]} \sup_{k\in\mathbb{N}} (x_{k+1,i} - x_i)^2  < + \infty$, where $\bm{x}:= (x_i) \in \mathbb{R}^d$ and $(\bm{x}_k)_{k\in\mathbb{N}} := ((x_{k,i}))_{k\in\mathbb{N}}$ is the sequence generated by Algorithm \ref{algo:1}.
\end{enumerate}
}
\end{assum}

The previous results in \cite[p.29]{chen2019}, \cite[p.18]{ada}, and \cite{iiduka2021} show that $(\mathsf{H}_k)_{k\in\mathbb{N}}$ in Table \ref{table:ex} satisfies (A4) and (A5). Assumption (A6) is assumed in 
\cite[p.1574]{nem2009}, \cite[Theorem 4.1]{adam}, \cite[p.2]{reddi2018}, 
and \cite[Theorem 2.1]{ada}.
If (A6) holds, then there exists a bounded set $X \subset \mathbb{R}^d$ such that $(\bm{x}_k)_{k\in\mathbb{N}} \subset X$.
Accordingly, the Lipschitz continuity of $\nabla f_i$, the nonemptiness of $S_i$, and (A6) imply that (G1) with $G := Ln \sqrt{dD}$ holds (see Proposition \ref{prop:1} in Appendix for details). 
We define 
\begin{align*}
h_0^* := \min_{i\in [d]} h_{0,i} \text{ and }
H := \max_{i\in [d]} H_i,
\end{align*}
where $h_{k,i}$ and $H_i$ are defined as in Assumption \ref{assum:1}.

\section{Main Results}\label{sec:main}
\subsection{Constant Learning Rate Rule}
The following theorem gives the relationship between batch size $s$ and the number of steps $K_\epsilon$ needed for (\ref{evaluation}) for Algorithm \ref{algo:1} 
with a constant learning rate $\alpha_k = \alpha$ (see Table \ref{constant} for the specific results in Theorem \ref{thm:1} with $G := Ln \sqrt{dD}$ (i.e., under condition (G1))).

\begin{thm}\label{thm:1}
Suppose that Assumptions \ref{assum:0} and \ref{assum:1} hold
and let $s, \epsilon > 0$.

{\em (i)}
Consider Algorithm \ref{algo:1} with
\begin{align*}
\alpha_k := \alpha \in (0,1] \text{ and } \beta_k := \beta \in [0,b] \subset [0,1).
\end{align*}
Then, for all $K \geq 1$,
\begin{align*}
\min_{k\in [K]} \mathbb{E}\left[ \| \nabla f(\bm{x}_k)  \|^2 \right]
\leq 
\underbrace{\frac{d D H}{2(1-b)\alpha}}_{A_\alpha} \frac{s}{K} 
+ 
\underbrace{\frac{G^2 \alpha}{2(1-b)(1-\gamma)^2 h_0^*}}_{B_\alpha} \frac{1}{s}
+
\underbrace{\frac{\sqrt{dD} G}{1-b} \beta}_{C_\beta}.
\end{align*}

{\em (ii)}
Consider Algorithm \ref{algo:1} with
\begin{align*}
\alpha_k := \alpha \in (0,1] \text{ and } \beta_k := \beta < \min \left\{ \frac{1-b}{\sqrt{dD}G} \epsilon^2, b    \right\}.
\end{align*} 
Then, the number of steps $K_\epsilon$ needed to achieve (\ref{evaluation}) is expressed as the following rational function of batch size $s$:
\begin{align}\label{function_c}
K_\epsilon (s) = \frac{A_\alpha s^2}{(\epsilon^2 - C_\beta)s - B_\alpha}
\quad\text{ } \left(s \in \left(\frac{B_\alpha}{\epsilon^2 - C_\beta}, + \infty \right) \right).
\end{align}
In particular, the minimum value of $K_\epsilon$ needed to achieve (\ref{evaluation}) is 
\begin{align*}
K_\epsilon(s^\star) 
=
\frac{4 A_\alpha B_\alpha}{(\epsilon^2 - C_\beta)^2}
= \frac{dDG^2 H}{(1-\gamma)^2 \{(1-b)\epsilon^2 - \sqrt{dD}G \beta \}^2 h_0^*}
\end{align*}
when 
\begin{align*}
s^\star 
=
\frac{2 B_\alpha}{\epsilon^2 - C_\beta}
= 
\frac{G^2 \alpha}{(1-\gamma)^2 \{(1-b)\epsilon^2 - \sqrt{dD}G\beta \}h_0^*}.
\end{align*}
\end{thm}

\subsubsection{Discussion of Theorem \ref{thm:1}}\label{sec:discussion_c}
Let us examine the results in Theorem \ref{thm:1} for SGD, N-Momentum, and Adam-type optimizers. 

\textbf{[Performance of Algorithm \ref{algo:1}]}
SGD is Algorithm \ref{algo:1} with $\beta = b = \gamma = 0$ and $h_0^* = H = 1$,
N-Momentum is Algorithm \ref{algo:1} with $\gamma = 0$ and $h_0^* = H = 1$,
and 
the Adam-type optimizer is Algorithm \ref{algo:1} with $\gamma \in [0,1)$ and $h_{k,i}$ defined by one of
$\sqrt{\hat{v}_{k,i}}$, $\sqrt{\tilde{v}_{k,i}}$, and $\sqrt{\hat{s}_{k,i}}$
(see Table \ref{table:ex}).
Theorem \ref{thm:1}(i) indicates that, for all $K \geq 1$, all $\alpha \in (0,1]$, all $\beta \in [0,b] \subset [0,1)$, and all $s > 0$,
\begin{align}\label{c_1}
\begin{split}
\min_{k\in [K]} \mathbb{E}\left[ \| \nabla f(\bm{x}_k)  \|^2 \right]
\leq
\begin{cases}
\frac{d D_{\mathrm{SGD}}}{2\alpha} \frac{s}{K} 
+ 
\frac{G^2 \alpha}{2} \frac{1}{s} &\text{ (SGD)},\\
\frac{d D_{\mathrm{NM}}}{2(1-b)\alpha} \frac{s}{K} 
+ 
\frac{G^2 \alpha}{2(1-b)} \frac{1}{s}
+
\frac{\sqrt{dD_{\mathrm{NM}}} G}{1-b} \beta &\text{ (N-Momentum)},\\
\frac{d D_{\mathrm{A}} H}{2(1-b)\alpha} \frac{s}{K} 
+ 
\frac{G^2 \alpha}{2(1-b)(1-\gamma)^2 h_0^*} \frac{1}{s}
+
\frac{\sqrt{dD_{\mathrm{A}}} G}{1-b} \beta &\text{ (Adam-type)}.
\end{cases}
\end{split}
\end{align}
Note that $D$ depends on the optimizer, which we distinguish by the notation 
$D_{\mathrm{SGD}}$, $D_{\mathrm{NM}}$, and $D_{\mathrm{A}}$.
For fixed $s$, if $\alpha$ and $\beta$ are sufficiently small, 
(\ref{c_1}) indicates that SGD, N-Momentum, and Adam-type optimizers have approximately $\mathcal{O}(1/K)$ convergence.
For fixed $s$ and $K$, if $\alpha$ is sufficiently small, the second term on the right-hand side of (\ref{c_1}) will be small, whereas the first term will be large.
Hence, there is no evidence that Algorithm \ref{algo:1} with a sufficiently small learning rate $\alpha$ would perform arbitrarily well. 
For fixed $\alpha$ and $K$,
if $s$ is sufficiently large, again the second term of the right-hand side of (\ref{c_1}) will be small and the first term will be large.
Hence, (\ref{c_1}) indicates that there is no evidence that Algorithm \ref{algo:1} with a large batch size $s$ performs better than with a smaller batch size.

\textbf{[Existence of optimal batch size]}
The function $K_\epsilon (s)$ defined by (\ref{function_c}) satisfies the following:
\begin{align*}
\frac{\mathrm{d} K_\epsilon (s)}{\mathrm{d}s}
\begin{cases}
< 0 \quad\text{ if } s \in \left(\frac{B_\alpha}{\epsilon^2 - C_\beta}, s^\star \right),\\ 
= 0 \quad\text{ if } s = s^\star = \frac{2B_\alpha}{\epsilon^2 - C_\beta},\\
> 0 \quad\text{ if } s \in (s^\star, + \infty).
\end{cases}
\end{align*}
The above shows that increasing the batch size initially decreases the number of steps $K_\epsilon$ needed to achieve (\ref{evaluation}). 
Then, there is an optimal batch size ($s = s^\star$) minimizing $K_\epsilon (s)$; thus increasing the batch size does not always decrease the number of steps $K_\epsilon$.

\textbf{[Comparison of optimal batch sizes]}
We assume that SGD, N-Momentum, and Adam-type optimizers all use the same $G$.
For example, under (G1), we have $G = Ln \sqrt{dD}$, where $D = \max \{D_{\mathrm{SGD}}, D_{\mathrm{NM}}, D_{\mathrm{A}}\}$.
From Theorem \ref{thm:1}(ii), we find that 
\begin{align}\label{c_2}
s_{\mathrm{C,SGD}}^\star 
= 
\frac{G^2 \alpha}{\epsilon^2}
\leq 
s_{\mathrm{C,NM}}^\star 
= 
\frac{G^2 \alpha}{(1-b)\epsilon^2 - \sqrt{dD_{\mathrm{NM}}} G\beta}.
\end{align}
This implies that N-Momentum exploits larger batches than SGD.
Moreover, if\footnote{$\gamma$ and $h_0^*$ can be chosen before implementing optimizers.
For example, let $\gamma = 0.9$, which is used in \cite{adam}.
Then, for all $i\in [d]$, we can set $h_{0,i} \leq 100$ (e.g., $h_0^* = h_{0,i} = 1$)
in order to satisfy (\ref{ams_c_3}).} 
\begin{align}\label{ams_c_3}
(1 - \gamma)^2 \leq \frac{1}{h_0^*},
\end{align} 
then we have that 
\begin{align}\label{c_2_1}
s_{\mathrm{C,SGD}}^\star 
= 
\frac{G^2 \alpha}{\epsilon^2}
\leq 
s_{\mathrm{C,A}}^\star 
= 
\frac{G^2 \alpha}{(1-\gamma)^2 \{(1-b)\epsilon^2 - \sqrt{dD_{\mathrm{A}}}G\beta \}h_0^*}.
\end{align}
Therefore, N-Momentum and Adam-type optimizers exploit larger batches than SGD.
Moreover, if (\ref{ams_c_3}) holds and if\footnote{We may assume that $D_{\mathrm{NM}} = D_{\mathrm{A}}$ in place of (\ref{Distance_0}).} 
\begin{align}\label{Distance_0}
D_{\mathrm{NM}} \leq D_{\mathrm{A}},
\end{align}
then 
\begin{align*}
s_{\mathrm{C,NM}}^\star 
\leq
s_{\mathrm{C,A}}^\star.
\end{align*}

\textbf{[Comparison of minimum numbers of steps]}
Theorem \ref{thm:1}(ii) guarantees that, if $\beta$ satisfies the condition in Theorem \ref{thm:1}(ii) and if
\begin{align}\label{b_1}
\beta \leq 
\frac{(1-b) \sqrt{D_{\mathrm{SGD}}} - \sqrt{D_{\mathrm{NM}}}}{\sqrt{d D_{\mathrm{SGD}} D_{\mathrm{NM}}}G}\epsilon^2,
\end{align}
then
\begin{align}\label{c_5}
K_\epsilon (s_{\mathrm{C,SGD}}^\star) 
= 
\frac{d D_{\mathrm{SGD}} G^2}{\epsilon^4}
\geq 
K_\epsilon (s_{\mathrm{C,NM}}^\star) 
= 
\frac{d D_{\mathrm{NM}} G^2}{\{(1-b)\epsilon^2 - \sqrt{dD_{\mathrm{NM}}}G\beta\}^2}.
\end{align}
Moreover, if $\beta$ satisfies the condition in Theorem \ref{thm:1}(ii) and if
\begin{align}\label{b_2}
\beta \leq 
\frac{(1-b)(1-\gamma) \sqrt{D_{\mathrm{SGD}} h_0^*} - \sqrt{D_{\mathrm{A}}H}}{(1-\gamma) \sqrt{d D_{\mathrm{SGD}} D_{\mathrm{A}} h_0^*}G}\epsilon^2,
\end{align}
then
\begin{align}\label{c_6}
K_\epsilon(s_{\mathrm{C,SGD}}^\star) 
\geq
K_\epsilon(s_{\mathrm{C,A}}^\star) 
= 
\frac{dD_{\mathrm{A}} G^2 H}{(1-\gamma)^2 \{(1-b)\epsilon^2 - \sqrt{dD_{\mathrm{A}}}G\beta \}^2 h_0^*}.
\end{align}
Additionally, if $\beta$ satisfies the condition in Theorem \ref{thm:1}(ii) and if 
\begin{align}\label{b_3}
\beta \leq 
\frac{(1-b)\{ (1-\gamma) \sqrt{D_{\mathrm{NM}} h_0^*} - \sqrt{D_{\mathrm{A}}H}\}}{\{(1-\gamma) \sqrt{d D_{\mathrm{NM}} D_{\mathrm{A}} h_0^*} - \sqrt{d D_{\mathrm{NM}} D_{\mathrm{A}} H} \} G}\epsilon^2,
\end{align}
then
\begin{align}\label{c_7}
K_\epsilon(s_{\mathrm{C,NM}}^\star) \geq K_\epsilon(s_{\mathrm{C,A}}^\star).
\end{align}
See (\ref{specific}) for more specific $\beta$ satisfying (\ref{b_3}).

\subsection{Diminishing Learning Rate Rule}\label{subsec:3.2}
The following theorem gives the relationships between batch size $s$ and the number of steps $K_\epsilon$ needed for (\ref{evaluation}) for Algorithm \ref{algo:1} 
with a diminishing learning rate $\alpha_k := \alpha/\sqrt{k}$
(see Table \ref{diminishing} for the specific results in Theorem \ref{thm:2} with $G := Ln \sqrt{dD}$
(i.e., under condition (G1))
and Theorem \ref{thm:3} for other results of Algorithm \ref{algo:1} with diminishing learning rates).

\begin{thm}\label{thm:2}
Suppose that Assumptions \ref{assum:0} and \ref{assum:1} hold
and also $s, \epsilon > 0$ and $\alpha \in (0,1]$.

{\em (i)}
Consider Algorithm \ref{algo:1} with
\begin{align*}
\alpha_k := \frac{\alpha}{\sqrt{k}} \text{ and } \beta_k := \beta \in [0,b] \subset [0,1).
\end{align*}
Then, for all $K \geq 1$,
\begin{align*}
\min_{k\in [K]} \mathbb{E}\left[ \| \nabla f(\bm{x}_k)  \|^2 \right]
\leq 
\underbrace{\frac{d D H}{2(1-b)\alpha}}_{A_\alpha} \frac{s}{\sqrt{K}} 
+ 
\underbrace{\frac{G^2 \alpha}{(1-b)(1-\gamma)^2 h_0^*}}_{B_\alpha} \frac{1}{s \sqrt{K}}
+
\underbrace{\frac{\sqrt{dD} G}{1-b} \beta}_{C_\beta}.
\end{align*}

{\em (ii)}
Consider Algorithm \ref{algo:1} with
\begin{align*}
\alpha_k := \frac{\alpha}{\sqrt{k}} \text{ and } \beta_k := \beta < \min \left\{ \frac{1-b}{\sqrt{dD}G} \epsilon^2, b    \right\}.
\end{align*} 
Then, the number of steps $K_\epsilon$ needed to achieve (\ref{evaluation}) is expressed as the following rational function of batch size $s$:
\begin{align}\label{function_d}
K_\epsilon (s) = \left\{ \frac{A_\alpha s^2 + B_\alpha}{(\epsilon^2 - C_\beta) s} \right\}^2.
\end{align}
In particular, the minimum value of $K_\epsilon$ needed to achieve (\ref{evaluation}) is 
\begin{align*}
K_\epsilon (s^\star) 
=
\frac{4 A_\alpha B_\alpha}{(\epsilon^2 - C_\beta)^2}
=
\frac{2 d D G^2 H}{(1-\gamma)^2\{ (1-b)\epsilon^2 - \sqrt{dD}G \beta \}^2 h_0^*}
\end{align*}
when 
\begin{align*}
s^\star 
=
\sqrt{\frac{B_\alpha}{A_\alpha}}
= 
\frac{\sqrt{2}G \alpha}{(1-\gamma) \sqrt{dDH h_0^*}}.
\end{align*}
\end{thm}

\subsubsection{Discussion of Theorem \ref{thm:2}}
Let us discuss the results in Theorem \ref{thm:2} and compare them with those in Theorem \ref{thm:1}  
for SGD, N-Momentum, and Adam-type optimizers. 

\textbf{[Performance of Algorithm \ref{algo:1}]}
Theorem \ref{thm:2}(i) indicates that
Algorithm \ref{algo:1} satisfies that, 
for all $K \geq 1$, all $\alpha \in (0,1]$, all $\beta \in [0,b]$, and all $s > 0$,
\begin{align*}
\min_{k\in [K]} \mathbb{E}\left[ \| \nabla f(\bm{x}_k)  \|^2 \right]
\leq
\begin{cases} 
\frac{d D_{\mathrm{SGD}}}{2\alpha} \frac{s}{\sqrt{K}}
+ 
\frac{G^2 \alpha}{s \sqrt{K}} &\text{ (SGD)},\\
\frac{d D_{\mathrm{NM}}}{2(1-b)\alpha} \frac{s}{\sqrt{K}}
+ 
\frac{G^2 \alpha}{(1-b)} \frac{1}{s \sqrt{K}}
+
\frac{\sqrt{dD_{\mathrm{NM}}} G}{1-b} \beta &\text{ (N-Momentum)},\\
\frac{d D_{\mathrm{A}} H}{2(1-b)\alpha} \frac{s}{\sqrt{K}}
+ 
\frac{G^2 \alpha}{(1-b)(1-\gamma)^2 h_0^*} \frac{1}{s \sqrt{K}}
+
\frac{\sqrt{dD_{\mathrm{A}}} G}{1-b} \beta &\text{ (Adam-type)}.
\end{cases}
\end{align*}
By a similar argument to that in Section \ref{sec:discussion_c},
SGD, N-Momentum, and Adam-type optimizers have approximately $\mathcal{O}(1/\sqrt{K})$ convergence (see also Theorem \ref{thm:3}, which indicates that Algorithm \ref{algo:1} with $\alpha_k = \alpha/\sqrt{k}$ and $\beta_k = \beta^k$ has an only $\mathcal{O}(1/\sqrt{K})$ convergence rate)
and that there is no evidence that Algorithm \ref{algo:1} with a large batch size $s$ performs better than with a smaller batch size.

\textbf{[Existence of optimal batch size]}
$K_\epsilon$ defined by (\ref{function_d}) guarantees that there exists $s^\star$ such that $\mathrm{d} K_\epsilon (s^\star)/\mathrm{d}s = 0$, the same as seen in Section \ref{sec:discussion_c} for Theorem 3.1. 
This implies that there is an optimal batch size ($s = s^\star$) such that
$K_\epsilon(s)$ is minimized, i.e., that increasing the batch size does not always decrease the number of steps $K_\epsilon$.

\textbf{[Comparison of optimal batch sizes]}
For simplicity, let us consider the case where (G1) holds.
Theorem \ref{thm:2}(ii) with $G = Ln \sqrt{dD}$ ensures that the optimal batch sizes for SGD, N-Momentum, and Adam-type optimizers with $\alpha_k = \alpha/\sqrt{k}$ and $\beta_k = \beta$ satisfy that
\begin{align*}
&s_{\mathrm{D,SGD}}^\star  
= 
\frac{\sqrt{2}G\alpha}{\sqrt{dD_{\mathrm{SGD}}}}
= \sqrt{2}Ln\alpha  
= 
\frac{\sqrt{2}G\alpha}{\sqrt{dD_{\mathrm{NM}}}}
=
s_{\mathrm{D,NM}}^\star.
\end{align*}
Furthermore, if\footnote{The definitions of $H$ and $h_0^*$ imply that $(1-\gamma) \sqrt{Hh_0^*} \leq (1-\gamma)H$. The condition $(1-\gamma)H \leq 1$ is sufficient to guarantee (\ref{H}).}
\begin{align}\label{H}
(1-\gamma)^2 \leq \frac{1}{H h_0^*},
\end{align}
then
\begin{align*}
&s_{\mathrm{D,SGD}}^\star  
=  
s_{\mathrm{D,NM}}^\star 
\leq
s_{\mathrm{D,A}}^\star
= \frac{\sqrt{2}Ln\alpha}{(1-\gamma) \sqrt{H h_0^*}}.
\end{align*}
Therefore, N-Momentum and Adam-type optimizers exploit the same sized or larger  batches than SGD. 
Here, we notice that 
$s_{\mathrm{C,SGD}}^\star$, $s_{\mathrm{C,NM}}^\star$, and $s_{\mathrm{C,A}}^\star$ defined as in (\ref{c_2}) and (\ref{c_2_1}) depend on $\epsilon$ and $\beta$, 
while 
$s_{\mathrm{D,SGD}}^\star$, $s_{\mathrm{D,NM}}^\star$, and 
$s_{\mathrm{D,A}}^\star$ do not depend on $\epsilon$  and $\beta$.

\textbf{[Comparison of minimum numbers of steps]}
Again, by a similar argument to that in Section \ref{sec:discussion_c}, 
the restrictions on $\beta$ (\ref{b_1}), (\ref{b_2}), and (\ref{b_3}) imply that (\ref{c_5}),  (\ref{c_6}), and (\ref{c_7}) hold, respectively, i.e., that
\begin{align*}
K_\epsilon (s_{\mathrm{D,A}}^\star) \leq K_\epsilon (s_{\mathrm{D,NM}}^\star)
\leq K_\epsilon (s_{\mathrm{D,SGD}}^\star).
\end{align*}
The previous studies \cite{adam,reddi2018,luo2019} used $\beta = 0.9$ or $0.99$, which is close to $1$, for adaptive methods. 
Meanwhile, a sufficient condition for $K_\epsilon (s_{\mathrm{D,A}}^\star) \leq K_\epsilon (s_{\mathrm{D,NM}}^\star)$ is (\ref{b_3}) with $D = D_{\mathrm{NM}} = D_{\mathrm{A}}$ and 
$G = Ln \sqrt{dD}$, i.e., 
\begin{align}\label{specific}
\beta \leq 
\frac{(1-b)\{ (1-\gamma) \sqrt{D_{\mathrm{NM}} h_0^*} - \sqrt{D_{\mathrm{A}}H}\}}{\{(1-\gamma) \sqrt{d D_{\mathrm{NM}} D_{\mathrm{A}} h_0^*} - \sqrt{d D_{\mathrm{NM}} D_{\mathrm{A}} H} \} G}\epsilon^2
=
\frac{(1-b)\epsilon^2}{\sqrt{dD}G}
= 
\frac{(1-b)\epsilon^2}{Ln dD},
\end{align}
which implies that adaptive methods using the above $\beta$ (which is small when the number of samples $n$ and the number of dimension $d$ are both large and the precision accuracy $\epsilon$ is small) are good for training deep neural networks in the sense that $K_\epsilon (s_{\mathrm{D,A}}^\star) \leq K_\epsilon (s_{\mathrm{D,NM}}^\star)$.

\section{Conclusion and Future Work}
The main contribution of this paper was to show that the number of steps $K_\epsilon (s)$ needed for nonconvex optimization, $\min_{k\in [K]} \mathbb{E}[\|\nabla f (\bm{x}_k)\|^2] \leq \epsilon^2$, of a deep learning optimizer is a rational function of batch size.
We showed that there exists an optimal batch size $s^\star$ such that $K_\epsilon(s)$ is minimized. 
This means that the optimizer using the optimal batch size $s^\star$ converges to a local minimizer of the sum of loss functions 
in at most $K_\epsilon(s^\star)$ steps and is most desirable for training deep neural networks. 
Hence, there is no guarantee that the optimizer with a sufficiently large batch size $s$ ($> s^\star$) would perform better than with a smaller batch size.
We also showed that the optimal batch size depends on the optimizer. 
In particular, it was shown that momentum and adaptive methods can exploit larger optimal batches
than can SGD and that, if we can set an appropriate momentum coefficient $\beta$, then 
momentum and adaptive methods reduce $K_\epsilon (s^\star)$ more than can SGD.

The results in this paper support theoretically the detailed numerical validations in recent papers \cite{shallue2019,zhang2019}.
The learning rate used in \cite[p.15]{shallue2019} decayed linearly,
which is distinctly  different from both constant and diminishing learning rates.   
In the future, we should check numerically the existences of optimal batch sizes of optimizers with not only constant but also diminishing learning rates to fully support all of the results in this paper.


\subsubsection*{Acknowledgments}
This work was supported by Japan Society for the Promotion of Science (JSPS) KAKENHI Grant Number 21K11773.

\bibliographystyle{spmpsci}      

\bibliography{iclr2022_conference}

\appendix
\section{Appendix}\label{appendix}
Unless stated otherwise, all relations between random variables
are supported to hold almost surely.
Let $S \in \mathbb{S}_{++}^d$. The $S$-inner product of $\mathbb{R}^d$ is defined for all $\bm{x}, \bm{y} \in \mathbb{R}^d$ by $\langle \bm{x},\bm{y} \rangle_S := \langle \bm{x}, S \bm{y} \rangle$ and the $S$-norm is defined by $\|\bm{x}\|_S := \sqrt{\langle \bm{x}, S \bm{x} \rangle}$. 
The history of process ${\xi}_0,{\xi}_1,\ldots$ to time step $k$ is denoted by ${\xi}_{[k]} = ({\xi}_0,{\xi}_1,\ldots,{\xi}_k)$.

\subsection{Sufficient Conditions for Assumption (A3)}
\begin{prop}\label{prop:1}
Assumption (A3) holds if each of the following holds:
\begin{enumerate}
\item[(G1)]
$X \subset \mathbb{R}^d$ is bounded, the gradient $\nabla f_i$ is Lipschitz continuous with Lipschitz constant $L_i$,
$S_i := \{ \bm{x}^* \in \mathbb{R}^d \colon \nabla f_i(\bm{x}^*) = \bm{0} \} \neq \emptyset$ ($i\in [n]$), 
where $L := \max_{i\in [n]} L_i$.
(If we define $G_{k,L} := \sup_{\bm{x} \in X}  \sum_{i\in \mathcal{S}_k} \|\nabla f_i(\bm{x}) \|$, then we can take $G := \sup_{k\in \mathbb{N}} G_{k,L}$.)
\item[(G2)]
$X \subset \mathbb{R}^d$ is bounded and closed. (If we define $G_{k} := \sup_{\bm{x} \in X}  \sum_{i\in \mathcal{S}_k} \|\nabla f_i(\bm{x}) \|$, then we can take $G := \sup_{k\in \mathbb{N}} G_{k}$.)\end{enumerate}
Under (A6), $G$ in (G1) and (G2) are respectively $G = Ln \sqrt{dD}$ and $G = n \tilde{G}$,
where $\tilde{G} := \max_{i\in [n]} \sup_{\bm{x} \in X} \|\nabla f_i (\bm{x})\|$.
\end{prop}

{\em Proof:}
The definition of $\nabla f_{\mathcal{S}_k}$ and the triangle inequality imply that, for all $\bm{x} \in \mathbb{R}^d$ and all $k\in\mathbb{N}$,
\begin{align}\label{ineq:1_0}
\left\|\nabla f_{\mathcal{S}_k}(\bm{x}) \right\|^2 
= 
\left\| \frac{1}{s} \sum_{i\in \mathcal{S}_k} \nabla f_i(\bm{x})\right\|^2
\leq 
\frac{1}{s^2} \left(\sum_{i\in \mathcal{S}_k} \left\| \nabla f_i(\bm{x})\right\| \right)^2.
\end{align} 
Suppose that (G1) holds.
Let $\bm{x}^* \in S_i$ ($i\in [n]$).
The Cauchy--Schwarz inequality and the Lipschitz continuity of $\nabla f_i$,
together with the definition of $L$, ensure that, for all $\bm{x} \in \mathbb{R}^d$ and all $i\in [n]$,
\begin{align*}
\|\nabla f_i(\bm{x})\| \leq \|\nabla f_i(\bm{x}^*)\| + L_i \|\bm{x} - \bm{x}^*\|
\leq L \|\bm{x} - \bm{x}^*\|.
\end{align*}
Accordingly, we have that, for all $\bm{x}\in X$ and all $k\in\mathbb{N}$, 
\begin{align*}
\sum_{i\in \mathcal{S}_k} \|\nabla f_i(\bm{x})\| 
\leq 
Ls \|\bm{x} - \bm{x}^*\|
\leq 
Ln \|\bm{x} - \bm{x}^*\|.
\end{align*}
Hence, $G_{k,L} \leq Ln \sup_{\bm{x}\in X} \|\bm{x} - \bm{x}^*\| < + \infty$.
Taking the expectation of (\ref{ineq:1_0}) thus implies (A3). 
Assumption (A6) implies that there exists a bounded set $X \subset \mathbb{R}^d$ such that 
$(\bm{x}_k)_{k\in\mathbb{N}} \subset X$. 
From 
$\|\bm{x}_k - \bm{x}^*\|^2 = \sum_{i\in [d]} (x_{k,i} - x_i)^2  \leq d D$,
we have that, for all $k\in\mathbb{N}$, 
\begin{align*}
G_{k,L} \leq Ln \sqrt{dD} =: G.
\end{align*} 

Suppose that (G2) holds. Since $\nabla f_i$ is continuous and $X$ is compact, 
we have that $G = \sup_{k\in\mathbb{N}} G_k < + \infty$.
Taking the expectation of (\ref{ineq:1_0}) thus implies (A3).
Assumption (A6) ensures that there exists a bounded, closed set $X \subset \mathbb{R}^d$ such that 
$(\bm{x}_k)_{k\in\mathbb{N}} \subset X$.
Define $\tilde{G}_i := \sup_{\bm{x} \in X} \|\nabla f_i(\bm{x})\| < + \infty$
and $\tilde{G} := \max_{i\in [n]} \tilde{G}_i$. 
Then, we have that, for all $\bm{x} \in X$, 
\begin{align*}
\sum_{i \in \mathcal{S}_k} \|\nabla f_i (\bm{x})\| \leq s \tilde{G} \leq n \tilde{G} =: G.
\end{align*}
This completes the proof.
$\Box$

\subsection{Examples of Algorithm \ref{algo:1}}
We list some examples of $\mathsf{H}_k \in \mathbb{S}_{++}^d \cap \mathbb{D}^d$ (step 5) in Algorithm \ref{algo:1}.

\begin{table}[htbp]
\caption{Examples of $\mathsf{H}_k \in \mathbb{S}_{++}^d \cap \mathbb{D}^d$ (step 5) in Algorithm \ref{algo:1} ($\delta,\zeta \in [0,1)$)}
\label{table:ex}
\centering
\begin{tabular}{l|l}
\hline
& $\mathsf{H}_k$ \\ \hline \hline
SGD   
& $\mathsf{H}_k$ is the identity matrix. \\
($\beta_k = \gamma = 0$)
&
\\ \hline
N-Momentum \cite{nest1983} 
&
$\mathsf{H}_k$ is the identity matrix. \\
($\gamma = 0$)
\\ \hline
AMSGrad \cite{chen2019} 
&
$\bm{v}_k = \delta \bm{v}_{k-1} + (1-\delta) \nabla f_{\mathcal{S}_k}(\bm{x}_k) \odot \nabla f_{\mathcal{S}_k}(\bm{x}_k)$ \\
($\gamma = 0$) &
$\hat{\bm{v}}_k = (\max \{ \hat{v}_{k-1,i}, v_{k,i} \})_{i=1}^d$ \\
&
$\mathsf{H}_k = \mathsf{diag} (\sqrt{\hat{v}_{k,i}})$ \\ \hline
AMSBound \cite{luo2019} & 
$\bm{v}_k = \delta \bm{v}_{k-1} + (1- \delta) \nabla f_{\mathcal{S}_k}(\bm{x}_k) \odot \nabla f_{\mathcal{S}_k}(\bm{x}_k)$ \\
($\gamma = 0$) &
$\hat{\bm{v}}_k = (\max \{ \hat{v}_{k-1,i}, v_{k,i} \})_{i=1}^d$ \\
&
$\tilde{\bm{v}}_k 
= \left(\mathrm{Clip}\left( \frac{1}{\sqrt{\hat{v}_{k,i}}}, l_k, u_k \right)^{-1} \right)_{i=1}^d$ \\
&
$\mathsf{H}_k = \mathsf{diag} (\sqrt{\tilde{v}_{k,i}})$ \\ \hline
AdaBelief \cite{ada} 
& 
$\tilde{\bm{s}}_k = (\nabla f_{\mathcal{S}_k}(\bm{x}_k) - \bm{m}_k) \odot 
(\nabla f_{\mathcal{S}_k}(\bm{x}_k) - \bm{m}_k)$ \\
($s_{k,i} \leq s_{k+1,i}$ is needed) &
$\bm{s}_k = \delta \bm{v}_{k-1} + (1-\delta) \tilde{\bm{s}}_k$ \\
&
$\hat{\bm{s}}_k = \frac{\bm{s}_k}{1- \zeta^k}$ \\
&
$\mathsf{H}_k = \mathsf{diag} (\sqrt{\hat{s}_{k,i}})$ \\ \hline
\end{tabular}
\end{table}

We define $\bm{x} \odot \bm{x}$ for $\bm{x} := (x_i)_{i=1}^d \in \mathbb{R}^d$ by $\bm{x} \odot \bm{x} := (x_i^2)_{i=1}^d \in \mathbb{R}^d$.
$\mathrm{Clip}(\cdot, l,u) \colon \mathbb{R} \to \mathbb{R}$ in AMSBound
($l,u \in \mathbb{R}$ with $l \leq u$ are given) is defined for all $x \in \mathbb{R}$ by
\begin{align*}
\mathrm{Clip}(x,l,u) := 
\begin{cases}
l &\text{ if } x < l,\\
x &\text{ if } l \leq x \leq u,\\
u &\text{ if } x > u.
\end{cases}
\end{align*}

\subsection{Lemmas and Theorem}
The following are the key lemmas to prove the main theorems in this paper.

\begin{lem}\label{lem:1}
Suppose that (A1) and (A2) hold and consider Algorithm \ref{algo:1}. Then, for all $\bm{x} \in \mathbb{R}^d$ and all $k\in\mathbb{N}$,
\begin{align*}
\mathbb{E}\left[\left\| \bm{x}_{k+1} - \bm{x} \right\|_{\mathsf{H}_k}^2\right]
&\leq
\mathbb{E}\left[\left\| \bm{x}_{k} - \bm{x} \right\|_{\mathsf{H}_k}^2\right]
+ \alpha_k^2 \mathbb{E} \left[ \left\|\bm{\mathsf{d}}_k \right\|_{\mathsf{H}_k}^2 \right]\\
&\quad + 2 \alpha_k 
 \left\{
\frac{\tilde{\beta}_k}{s \tilde{\gamma}_k} 
 {\mathbb{E} \left[\left\langle \bm{x} - \bm{x}_k, \nabla f (\bm{x}_k) \right\rangle \right]}
+ \frac{\beta_k}{\tilde{\gamma}_k} \mathbb{E} \left[ \left\langle \bm{x} - \bm{x}_k,\bm{m}_{k-1} \right\rangle \right]
 \right\}, 
\end{align*}
where $\tilde{\beta}_k := 1 - \beta_k$ and $\tilde{\gamma}_k := 1 - \gamma^{k+1}$.
\end{lem}

{\em Proof:} 
Let $\bm{x}\in \mathbb{R}^d$ and $k\in \mathbb{N}$.
The definition of $\bm{x}_{k+1}$ implies that
\begin{align*}
\left\|\bm{x}_{k+1} - \bm{x} \right\|_{\mathsf{H}_k}^2
\leq
\left\| \bm{x}_k - \bm{x} \right\|_{\mathsf{H}_k}^2 
+ 2 \alpha_k \left\langle \bm{x}_k - \bm{x}, \bm{\mathsf{d}}_k \right\rangle_{\mathsf{H}_k}
+ \alpha_k^2 \left\|\bm{\mathsf{d}}_k \right\|_{\mathsf{H}_k}^2. 
\end{align*}
Moreover, the definitions of $\bm{\mathsf{d}}_k$, $\bm{m}_k$, and $\hat{\bm{m}}_k$ ensure that
\begin{align*}
\left\langle \bm{x}_k - \bm{x}, \bm{\mathsf{d}}_k \right\rangle_{\mathsf{H}_k}=
\frac{1}{{\tilde{\gamma}}_k}\left\langle \bm{x} - \bm{x}_k, \bm{m}_k \right\rangle
=
\frac{\beta_k}{{\tilde{\gamma}}_k} \left\langle \bm{x} - \bm{x}_k, \bm{m}_{k-1} \right\rangle 
+
\frac{\tilde{\beta}_k}{{\tilde{\gamma}}_k} \left\langle \bm{x} - \bm{x}_k, \nabla f_{\mathcal{S}_k}(\bm{x}_k) \right\rangle,
\end{align*}
where $\tilde{\beta}_k := 1 - \beta_k$ and ${\tilde{\gamma}}_k := 1 - {\gamma}^{k+1}$. 
Hence, 
\begin{align}\label{ineq:004}
\begin{split}
\left\|\bm{x}_{k+1} - \bm{x} \right\|_{\mathsf{H}_k}^2
&\leq
\left\| \bm{x}_k -\bm{x} \right\|_{\mathsf{H}_k}^2
+ 
2 \alpha_k \left\{
\frac{\beta_k}{{\tilde{\gamma}}_k} \left\langle \bm{x} - \bm{x}_k, \bm{m}_{k-1} \right\rangle
+ \frac{\tilde{\beta}_k}{{\tilde{\gamma}}_k} \left\langle \bm{x} - \bm{x}_k, \nabla f_{\mathcal{S}_k} (\bm{x}_k) \right\rangle
\right\}\\
&\quad + \alpha_k^2 \left\| \bm{\mathsf{d}}_k \right\|_{\mathsf{H}_k}^2.
\end{split}
\end{align}
Meanwhile, the relationship between the expectation of the stochastic gradient vector $\nabla f_{\mathcal{S}_k}(\bm{x})$ and the full gradient vector $\nabla f(\bm{x})$ is as follows:
For all $\bm{x} \in \mathbb{R}^d$, 
\begin{align}\label{equation}
\begin{split}
\mathbb{E}\left[ \nabla f_{\mathcal{S}_k}(\bm{x}) \right]
=
\mathbb{E}\left[ \frac{1}{s} \sum_{i\in \mathcal{S}_k} \nabla f_i (\bm{x}) \right]
= 
\frac{1}{s} \mathbb{E}\left[ \nabla f_i (\bm{x}) \right]
= 
\frac{1}{s} \nabla f (\bm{x}),
\end{split}
\end{align}
where the first equation comes from (A2), the second equation comes from the existence of $T$ such that $[n] = \cup_{k=1}^T \mathcal{S}_k$, and the third equation comes from (A1).
Condition (\ref{equation}) guarantees that
\begin{align*}
\mathbb{E} \left[\left\langle \bm{x} - \bm{x}_k, \nabla f_{\mathcal{S}_k} (\bm{x}_k) \right\rangle \right]
&=
\mathbb{E} 
 \left[ 
 \mathbb{E} \left[\left\langle \bm{x} - \bm{x}_k, \nabla f_{\mathcal{S}_k} (\bm{x}_k) \right\rangle | {\xi}_{[k-1]} \right]
 \right]\\
&= 
\mathbb{E} 
 \left[ 
 \left\langle \bm{x} - \bm{x}_k, \mathbb{E} \left[ \nabla f_{\mathcal{S}_k} (\bm{x}_k) | {\xi}_{[k-1]} \right] \right\rangle
 \right]\\
&=
\frac{1}{s} \mathbb{E} \left[\left\langle \bm{x} - \bm{x}_k, {\nabla f} (\bm{x}_k) \right\rangle \right].
\end{align*}
Therefore, the lemma follows by taking the expectation of (\ref{ineq:004}).
$\Box$

\begin{lem}\label{lem:bdd}
Algorithm \ref{algo:1} satisfies that, under (A3),  for all $k\in\mathbb{N}$,
\begin{align*}
\mathbb{E}\left[ \|\bm{m}_k\|^2 \right] 
\leq 
\frac{G^2}{s^2}.
\end{align*}
Under (A3) and (A4), for all $k\in\mathbb{N}$,
\begin{align*}
\mathbb{E}\left[ \|\bm{\mathsf{d}}_k\|_{\mathsf{H}_k}^2 \right] 
\leq 
\frac{G^2}{(1-\gamma)^2 h_0^* s^2},
\end{align*}
where $h_0^* := \min_{i\in [d]} h_{0,i}$.
\end{lem}

{\em Proof:}
The convexity of $\|\cdot\|^2$, together with the definition 
of $\bm{m}_k$ and (A3), guarantees that , for all $k\in\mathbb{N}$,
\begin{align*}
\mathbb{E}\left[ \|\bm{m}_k\|^2 \right]
&\leq \beta_k \mathbb{E}\left[ \|\bm{m}_{k-1} \|^2 \right] + (1-\beta_k)\mathbb{E}\left[ \|\nabla f_{\mathcal{S}_k} (\bm{x}_k) \|^2 \right]\\
&\leq 
\beta_k \mathbb{E} \left[ \|\bm{m}_{k-1} \|^2 \right] + (1-\beta_k) \frac{G^2}{s^2}.
\end{align*}
Induction thus ensures that , for all $k\in\mathbb{N}$,
\begin{align}\label{induction}
\mathbb{E} \left[\|\bm{m}_k \|^2 \right] \leq  
\max \left\{ \|\bm{m}_{-1}\|^2, \frac{G^2}{s^2} \right\} 
= \frac{G^2}{s^2},
\end{align}
where $\bm{m}_{-1} = \bm{0}$ is used.
For $k\in\mathbb{N}$, $\mathsf{H}_k \in \mathbb{S}_{++}^d$ guarantees the existence of a unique matrix $\overline{\mathsf{H}}_k \in \mathbb{S}_{++}^d$ such that $\mathsf{H}_k 
= \overline{\mathsf{H}}_k^2$ \cite[Theorem 7.2.6]{horn}. 
We have that, for all $\bm{x}\in\mathbb{R}^d$, 
$\|\bm{x}\|_{\mathsf{H}_k}^2 = \| \overline{\mathsf{H}}_k \bm{x} \|^2$.
Accordingly, the definitions of $\bm{\mathsf{d}}_k$ and $\hat{\bm{m}}_k$
imply that, for all $k\in\mathbb{N}$, 
\begin{align*}
\mathbb{E} \left[ \| \bm{\mathsf{d}}_k \|_{\mathsf{H}_k}^2 \right]
= 
\mathbb{E} \left[ \left\| \overline{\mathsf{H}}_k^{-1} \mathsf{H}_k\bm{\mathsf{d}}_k \right\|^2 \right]
\leq 
\frac{1}{{\tilde{\gamma}}_k^2} \mathbb{E} \left[ \left\| \overline{\mathsf{H}}_k^{-1} \right\|^2 \|\bm{m}_k \|^2 \right]
\leq 
\frac{1}{(1 - \gamma)^2} \mathbb{E} \left[ \left\| \overline{\mathsf{H}}_k^{-1} \right\|^2 \|\bm{m}_k \|^2 \right],
\end{align*}
where 
\begin{align*}
\left\| \overline{\mathsf{H}}_k^{-1} \right\| = \left\| \mathsf{diag}\left(h_{k,i}^{-\frac{1}{2}} \right) \right\| = {\max_{i\in [d]} h_{k,i}^{-\frac{1}{2}}} 
\end{align*}
and ${\tilde{\gamma}}_k := 1 - {\gamma}^{k+1} \geq 1 - {\gamma}$.
Moreover, (A4) ensures that, for all $k \in \mathbb{N}$, 
\begin{align*}
h_{k,i} \geq h_{0,i} \geq h_0^* := \min_{i\in [d]} h_{0,i}.
\end{align*}
Hence, (\ref{induction}) implies that , for all $k\in \mathbb{N}$,
\begin{align*}
\mathbb{E} \left[\| \bm{\mathsf{d}}_k \|_{\mathsf{H}_k}^2 \right] \leq 
\frac{G^2}{(1-{\gamma})^2 h_0^* s^2},
\end{align*}
completing the proof.
$\Box$

We are in the position to prove the following theorem, which leads to Theorems \ref{thm:1}, \ref{thm:2}, and \ref{thm:3}.

\begin{thm}\label{thm:main}
Suppose that Assumptions \ref{assum:0} and \ref{assum:1} hold and consider Algorithm \ref{algo:1}. Let $(\delta_k)_{k\in \mathbb{N}} \subset (0,+\infty)$ be the sequence defined by $\delta_k := \alpha_k \tilde{\beta_k}/\tilde{\gamma}_k$ and $V_k (\bm{x}) := \mathbb{E}[ \langle \bm{x}_k - \bm{x}, \nabla f (\bm{x}_k) \rangle ]$ for all $\bm{x} \in \mathbb{R}^d$ and all $k\in\mathbb{N}$. Assume that $(\delta_k)_{k\in\mathbb{N}}$ is monotone decreasing. 
Then, for all $\bm{x} \in \mathbb{R}^d$ and all $K\geq 1$,
\begin{align*}
\sum_{k=1}^K V_k(\bm{x})
\leq
\frac{ds DH}{2 \tilde{b} \alpha_K}
+
\frac{G^2}{2 \tilde{b} \tilde{\gamma}^2 h_0^* s} \sum_{k=1}^K \alpha_k
+
\frac{\sqrt{dD} G}{\tilde{b}} \sum_{k=1}^K \beta_k,
\end{align*}
where $\tilde{b} := 1 - b$, $\tilde{\gamma} := 1 - \gamma$, $D$ and $H_i$ are defined as in Assumption \ref{assum:1}, and $H := \max_{i\in [d]} H_i$.
\end{thm}

{\em Proof:} 
Let $\bm{x}\in \mathbb{R}^d$. Lemma \ref{lem:1} guarantees that, for all $k\in\mathbb{N}$,
\begin{align*}
V_k (\bm{x})
&\leq
\frac{s}{2\delta_k}
\left\{ 
\mathbb{E}\left[\left\| \bm{x}_{k} - \bm{x} \right\|_{\mathsf{H}_k}^2\right]
-
\mathbb{E}\left[\left\| \bm{x}_{k+1} - \bm{x} \right\|_{\mathsf{H}_k}^2\right]
\right\}
+ \frac{s \alpha_k {\tilde{\gamma}}_k}{2 \tilde{\beta}_k} \mathbb{E} \left[ \left\|\bm{\mathsf{d}}_k \right\|_{\mathsf{H}_k}^2 \right]\\
&\quad + \frac{s \beta_k}{\tilde{\beta_k}} \mathbb{E} \left[ \left\langle \bm{x} - \bm{x}_k, \bm{m}_{k-1} \right\rangle \right].
\end{align*}
Summing the above inequality from $k=1$ to $K \geq 1$ implies that 
\begin{align}\label{key}
\begin{split}
{\sum_{k=1}^K V_k (\bm{x})}
&\leq
\frac{1}{2} \underbrace{\sum_{k=1}^K \frac{s}{\delta_k}
\left\{ 
\mathbb{E}\left[\left\| \bm{x}_{k} - \bm{x} \right\|_{\mathsf{H}_k}^2\right]
-
\mathbb{E}\left[\left\| \bm{x}_{k+1} - \bm{x} \right\|_{\mathsf{H}_k}^2\right]
\right\}}_{\Delta_K}\\
&\quad 
+ \frac{1}{2} \underbrace{\sum_{k=1}^K \frac{s \alpha_k \tilde{\gamma}_k}{\tilde{\beta}_k} \mathbb{E} \left[ \left\|\bm{\mathsf{d}}_k \right\|_{\mathsf{H}_k}^2 \right]}_{A_K}
+ \underbrace{
\sum_{k=1}^K \frac{s \beta_k}{\tilde{\beta}_k} \mathbb{E} \left[ \left\langle \bm{x} - \bm{x}_k,\bm{m}_{k-1} \right\rangle \right]}_{B_K}.
\end{split}
\end{align} 
Let us define $\theta_k := \delta_k/s$. From the definition of $\Delta_K$ and $\mathbb{E} [ \| \bm{x}_{k+1} - \bm{x} \|_{\mathsf{H}_{k}}^2]/\theta_k \geq 0$, 
\begin{align}\label{LAM}
\Delta_K
\leq
\frac{\mathbb{E}\left[\left\| \bm{x}_{1} - \bm{x} \right\|_{\mathsf{H}_{1}}^2\right]}{\theta_1}
+
\underbrace{
\sum_{k=2}^K \left\{
\frac{\mathbb{E}\left[\left\| \bm{x}_{k} - \bm{x} \right\|_{\mathsf{H}_{k}}^2\right]}{\theta_k}
-
\frac{\mathbb{E}\left[\left\| \bm{x}_{k} - \bm{x} \right\|_{\mathsf{H}_{k-1}}^2\right]}{\theta_{k-1}} 
\right\}
}_{\tilde{\Delta}_K}.
\end{align}
Since $\overline{\mathsf{H}}_k \in \mathbb{S}_{++}^d$ exists such that $\mathsf{H}_k = \overline{\mathsf{H}}_k^2$, we have $\|\bm{x}\|_{\mathsf{H}_k}^2 = \| \overline{\mathsf{H}}_k \bm{x} \|^2$ for all $\bm{x}\in\mathbb{R}^d$. Accordingly, we have 
\begin{align*}
\tilde{\Delta}_K 
=
\mathbb{E} \left[ 
\sum_{k=2}^K 
\left\{
\frac{\left\| \overline{\mathsf{H}}_{k} (\bm{x}_{k} - \bm{x}) \right\|^2}{\theta_k}
-
\frac{\left\| \overline{\mathsf{H}}_{k-1} (\bm{x}_{k} - \bm{x}) \right\|^2}{\theta_{k-1}}
\right\}
\right].
\end{align*}
From $\overline{\mathsf{H}}_{k} = \mathsf{diag}(\sqrt{h_{k,i}})$,
we have that, for all $\bm{x} = (x_i)_{i=1}^d \in \mathbb{R}^d$, 
$\| \overline{\mathsf{H}}_{k} \bm{x} \|^2
=
\sum_{i=1}^d h_{k,i} x_i^2$. 
Hence, for all $K\geq 2$,
\begin{align*}
\tilde{\Delta}_K 
= 
\mathbb{E} \left[ 
\sum_{k=2}^K
\sum_{i=1}^d 
\left(
\frac{h_{k,i}}{\theta_k}
-
\frac{h_{k-1,i}}{\theta_{k-1}}
\right)
(x_{k,i} - x_i)^2
\right].
\end{align*}
Accordingly, from (A4) and the monotone decrease of $(\delta_k)_{k\in\mathbb{N}}$, we have that, for all $k \geq 1$ and all $i\in [d]$,
\begin{align*}
\frac{h_{k,i}}{\theta_k} - \frac{h_{k-1,i}}{\theta_{k-1}} \geq 0.
\end{align*} 
Moreover, from (A6), $D := \max_{i \in [d]} \sup_{k\in\mathbb{N}}
(x_{k,i} - x_i)^2 < + \infty$. Accordingly, for all $K \geq 2$,
\begin{align*}
\tilde{\Delta}_K
\leq
D
\mathbb{E} \left[ 
\sum_{k=2}^K
\sum_{i=1}^d 
\left(
\frac{h_{k,i}}{\theta_k}
-
\frac{h_{k-1,i}}{\theta_{k-1}}
\right)
\right]
= 
D
\mathbb{E} \left[ 
\sum_{i=1}^d
\left(
\frac{h_{K,i}}{\theta_K}
-
\frac{h_{1,i}}{\theta_{1}}
\right)
\right].
\end{align*}
Therefore, (\ref{LAM}), $\mathbb{E} [\| \bm{x}_{1} - \bm{x}\|_{\mathsf{H}_{1}}^2]/\theta_1 \leq D \mathbb{E} [ \sum_{i=1}^d h_{1,i}/\theta_1]$, and (A5) imply, for all $K\in \mathbb{N}$,
\begin{align*}
\Delta_K
\leq
D \mathbb{E} \left[ \sum_{i=1}^d \frac{h_{1,i}}{\theta_1} \right]
+
D
\mathbb{E} \left[
\sum_{i=1}^d 
\left(
\frac{h_{K,i}}{\theta_K}
-
\frac{h_{1,i}}{\theta_{1}}
\right)
\right]
=
\frac{D}{\theta_K}
\mathbb{E} \left[
\sum_{i=1}^d 
h_{K,i}
\right]
\leq 
\frac{D}{\theta_K}
\sum_{i=1}^d 
H_{i},
\end{align*}
which, together with $\theta_K := \alpha_K (1-\beta_K)/(s (1 - {\gamma}^{K+1})) \geq \tilde{b} \alpha_K / s$ and $H = \max_{i\in [d]} H_i$, implies 
\begin{align}\label{L} 
\frac{1}{2} \Delta_K 
\leq 
\frac{d s D H}{2 \tilde{b} \alpha_K}.
\end{align}
Lemma \ref{lem:bdd} implies that , for all $K\in\mathbb{N}$,
\begin{align*}
A_K := \sum_{k=1}^K \frac{s \alpha_k \tilde{\gamma}_k}{\tilde{\beta}_k} \mathbb{E} \left[ \left\|\bm{\mathsf{d}}_k \right\|_{\mathsf{H}_k}^2 \right] 
\leq 
\sum_{k=1}^K \frac{s \alpha_k \tilde{\gamma}_k}{\tilde{\beta}_k} \frac{G^2}{\tilde{\gamma}^2 h_0^* s^2},
\end{align*}
which, together with $\tilde{\gamma}_k \leq 1$ and $\beta_k \leq b$, implies that
\begin{align}\label{D}
\frac{1}{2} A_K \leq 
\frac{G^2}{2 \tilde{b} \tilde{\gamma}^2 h_0^* s} \sum_{k=1}^K \alpha_k.
\end{align}
Lemma \ref{lem:bdd} and Jensen's inequality ensure that, for all $k\in\mathbb{N}$,
\begin{align*}
\mathbb{E}\left[ \|\bm{m}_k\| \right] 
\leq 
\frac{G}{s}.
\end{align*}
The Cauchy-Schwarz inequality and (A6) guarantee that, for all $K\in\mathbb{N}$,
\begin{align}\label{B}
B_K 
:= 
\sum_{k=1}^K \frac{s \beta_k}{\tilde{\beta}_k} \mathbb{E} \left[ \left\langle \bm{x} - \bm{x}_k,\bm{m}_{k-1} \right\rangle \right]
\leq 
\sum_{k=1}^K \frac{s\sqrt{dD}\beta_k}{\tilde{b}} \mathbb{E} \left[
\left\|\bm{m}_{k-1} \right\|
\right]
\leq
\frac{\sqrt{dD}G}{\tilde{b}} \sum_{k=1}^K \beta_k.
\end{align}
Therefore, (\ref{key}), (\ref{L}), (\ref{D}), and (\ref{B}) lead to the assertion in Theorem \ref{thm:main}. This completes the proof.
$\Box$

\subsection{Proof of Theorem \ref{thm:1}}
(i) 
Theorem \ref{thm:main}, together with $\alpha_k = \alpha$ and $\beta_k = \beta$, guarantees that , for all $K \geq 1$ and all $\bm{x}\in \mathbb{R}^d$,
\begin{align}\label{ineq:1}
\frac{1}{K} \sum_{k=1}^K V_k(\bm{x})
\leq
\frac{d D H}{2 \tilde{b} \alpha}\frac{s}{K}
+
\frac{G^2 \alpha}{2 \tilde{b} \tilde{\gamma}^2 h_0^*} \frac{1}{s} 
+
\frac{\sqrt{dD} G}{\tilde{b}} \beta.
\end{align}
Moreover,
there exists $m \in [K]$ such that, for all $\bm{x} \in \mathbb{R}^d$,
\begin{align}\label{ineq:2}
\mathbb{E}\left[ \langle \bm{x}_m - \bm{x}, \nabla f (\bm{x}_m) \rangle \right]
= V_m (\bm{x})
= \min_{k\in [K]} V_k(\bm{x}) \leq \frac{1}{K} \sum_{k=1}^K V_k(\bm{x}).
\end{align}
Setting $\bm{x} = \bm{x}_m - \nabla f (\bm{x}_m)$, together with (\ref{ineq:1}) and (\ref{ineq:2}), guarantees that 
\begin{align}\label{ineq:3}
\min_{k\in [K]} \mathbb{E}\left[ \|\nabla f (\bm{x}_k)\|^2 \right]
\leq
\mathbb{E}\left[ \|\nabla f (\bm{x}_m)\|^2 \right]
\leq
\underbrace{\frac{d D H}{2 \tilde{b} \alpha}}_{A_\alpha} \frac{s}{K}
+
\underbrace{\frac{G^2 \alpha}{2 \tilde{b} \tilde{\gamma}^2 h_0^*}}_{B_\alpha} \frac{1}{s} 
+
\underbrace{\frac{\sqrt{dD} G}{\tilde{b}} \beta}_{C_\beta}.
\end{align} 

(ii)
A sufficient condition for (\ref{evaluation}), i.e., 
\begin{align*}
\min_{k\in [K]} \mathbb{E}\left[ \|\nabla f (\bm{x}_k)\|^2 \right]
\leq \epsilon^2
\end{align*} 
is that the right-hand side of (\ref{ineq:3}) is equal to $\epsilon^2$, i.e., 
\begin{align*}
A_\alpha s^2 + B_\alpha K + (C_\beta - \epsilon^2) sK = 0,
\end{align*}
which implies that 
\begin{align*}
K(s) = \frac{A_\alpha s^2}{(\epsilon^2 - C_\beta)s - B_\alpha}
\quad\text{ } \left(s \in \left(\frac{B_\alpha}{\epsilon^2 - C_\beta}, + \infty \right) \right),
\end{align*}
where $\epsilon^2 - C_\beta > 0$ is guaranteed from $\beta < \tilde{b}\epsilon^2/\sqrt{dD}G$.
We have that
\begin{align*}
\frac{\mathrm{d} K(s)}{\mathrm{d}s}
=
\frac{A_\alpha s}{\{(C_\beta - \epsilon^2)s + B_\alpha\}^2}
\left\{(\epsilon^2 - C_\beta)s - 2 B_\alpha \right\} 
\begin{cases}
\displaystyle{< 0 \quad\text{ if } s \in \left(\frac{B_\alpha}{\epsilon^2 - C_\beta},s^\star \right),}\\ 
\displaystyle{= 0 \quad\text{ if } s = s^\star = \frac{2 B_\alpha}{\epsilon^2 - C_\beta},}\\
\displaystyle{> 0 \quad\text{ if } s \in (s^\star, + \infty).}
\end{cases}
\end{align*}
Hence, $K(s)$ attains the minimum $K(s^\star)$ when $s = s^\star$.
$\Box$

\subsection{Proof of Theorem \ref{thm:2}}
(i) 
Theorem \ref{thm:main}, together with $\alpha_k = \alpha/\sqrt{k}$ and $\beta_k = \beta$, guarantees that, for all $K \geq 1$ and all $\bm{x}\in \mathbb{R}^d$,
\begin{align}\label{ineq:2_1}
\begin{split}
\frac{1}{K} \sum_{k=1}^K V_k(\bm{x})
&\leq
\frac{d D H}{2 \tilde{b}}\frac{s}{\alpha_K K}
+
\frac{G^2}{2 \tilde{b} \tilde{\gamma}^2 h_0^*}\frac{1}{s K} \sum_{k=1}^K \alpha_k
+
\frac{\sqrt{dD} G}{\tilde{b}} \beta\\
&\leq
\frac{d D H}{2 \tilde{b} \alpha}\frac{s}{\sqrt{K}}
+
\frac{G^2 \alpha}{\tilde{b} \tilde{\gamma}^2 h_0^*}\frac{1}{s \sqrt{K}}
+
\frac{\sqrt{dD} G}{\tilde{b}} \beta,
\end{split}
\end{align}
where we use 
\begin{align*}
\frac{1}{K} \sum_{k=1}^K \frac{1}{\sqrt{k}}
\leq 
\frac{1}{K} \left(1 + \int_1^K \frac{\mathrm{d} t}{\sqrt{t}}  \right)
= 
\frac{1}{K} \left(2 \sqrt{K} -1 \right)
\leq 
\frac{2}{\sqrt{K}}. 
\end{align*}
An argument similar to the one for showing (\ref{ineq:2}) and (\ref{ineq:3}) ensures that
(\ref{ineq:2_1}) implies that 
\begin{align}\label{ineq:2_3}
\min_{k\in [K]} \mathbb{E}\left[ \|\nabla f (\bm{x}_k)\|^2 \right]
\leq
\underbrace{\frac{d D H}{2 \tilde{b} \alpha }}_{A_\alpha} \frac{s}{\sqrt{K}}
+
\underbrace{\frac{G^2 \alpha}{\tilde{b} \tilde{\gamma}^2 h_0^*}}_{B_\alpha} \frac{1}{s \sqrt{K}} 
+
\underbrace{\frac{\sqrt{dD} G}{\tilde{b}} \beta}_{C_\beta}.
\end{align} 

(ii)
A sufficient condition for (\ref{evaluation}), i.e., 
\begin{align*}
\min_{k\in [K]} \mathbb{E}\left[ \|\nabla f (\bm{x}_k)\|^2 \right]
\leq \epsilon^2
\end{align*} 
is that the right-hand side of (\ref{ineq:2_3}) is equal to $\epsilon^2$, i.e., 
\begin{align*}
A_\alpha s^2 + (C_\beta - \epsilon^2) s \sqrt{K} + B_\alpha = 0,
\end{align*}
which implies that 
\begin{align*}
K(s) = \left\{ \frac{A_\alpha s^2 + B_\alpha}{(\epsilon^2 - C_\beta) s} \right\}^2,
\end{align*}
where $\epsilon^2 - C_\beta > 0$ is guaranteed from $\beta < \tilde{b}\epsilon^2/\sqrt{dD}G$.
We have that
\begin{align*}
\frac{\mathrm{d} K(s)}{\mathrm{d}s}
=
\frac{2 (A_\alpha s^2 +B_\alpha)}{(\epsilon^2 - C_\beta)^2 s^3}
(A_\alpha s^2  - B_\alpha) 
\begin{cases}
\displaystyle{< 0 \quad\text{ if } s \in ( 0, s^\star ),}\\ 
\displaystyle{= 0 \quad\text{ if } s = s^\star = \sqrt{\frac{B_\alpha}{A_\alpha}},}\\
\displaystyle{> 0 \quad\text{ if } s \in (s^\star, + \infty),}
\end{cases}
\end{align*}
which implies that $K(s)$ attains the minimum $K(s^\star)$ when $s = s^\star$.
$\Box$

\subsection{Relationship between $s$ and $K_\epsilon(s)$ for Algorithm \ref{algo:1} with diminishing learning rates}
The following is a result for Algorithm \ref{algo:1} with diminishing sequences $\alpha_k$ and $\beta_k$.

\begin{thm}\label{thm:3}
Suppose that Assumptions \ref{assum:0} and \ref{assum:1} hold
and let $s, \epsilon > 0$, $\alpha \in (0,1]$, and $\beta \in [0,b] \subset [0,1)$.

{\em (i)}
Consider Algorithm \ref{algo:1} with
\begin{align*}
\alpha_k := \frac{\alpha}{\sqrt{k}} \text{ and } \beta_k := \beta^k.
\end{align*}
Then, for all $K \geq 1$, 
\begin{align*}
\min_{k\in [K]} \mathbb{E}\left[ \| \nabla f(\bm{x}_k)  \|^2 \right]
\leq 
\underbrace{\frac{d D H}{2(1-b)\alpha}}_{A_\alpha} \frac{s}{\sqrt{K}} 
+ 
\underbrace{\frac{G^2 \alpha}{(1-b)(1-\gamma)^2 h_0^*}}_{B_\alpha} \frac{1}{s \sqrt{K}}
+
\underbrace{\frac{\beta \sqrt{dD} G}{(1-b) (1 - \beta)}}_{C_\beta} \frac{1}{K}.
\end{align*}

{\em (ii)}
The number of steps $K_\epsilon$ needed to achieve (\ref{evaluation}) is expressed as the following rational function of batch size $s$:
\begin{align*}
K_\epsilon (s) = \left\{ \frac{(A_\alpha s^2 + B_\alpha) + \sqrt{(A_\alpha s^2 + B_\alpha)^2 + 4 \epsilon^2 C_\beta s^2}}{2 \epsilon^2 s} \right\}^2.
\end{align*}
In particular, the minimum value of $K_\epsilon$ needed to achieve (\ref{evaluation}) is 
\begin{align*}
K_\epsilon (s^\star) 
&=
\left\{\frac{\sqrt{A_\alpha B_\alpha} + \sqrt{A_\alpha B_\alpha + \epsilon^2 C_\beta} }{\epsilon^2} \right\}^2\\
&=
\frac{\left\{ \sqrt{(1-\beta)dDH}G + \sqrt{\left((1-\beta)dDGH + 2 (1-b)\beta (1-\gamma)^2 \epsilon^2 \sqrt{dD}h_0^* \right)G} \right\}^2}{2 (1-b)^2 (1-\beta)(1-\gamma)^2 \epsilon^4 h_0^*}
\end{align*}
when 
\begin{align*}
s^\star 
=
\sqrt{\frac{B_\alpha}{A_\alpha}}
= 
\frac{G \alpha}{(1-\gamma) \sqrt{dDH h_0^*}}.
\end{align*}
\end{thm}

{\em Proof:}
(i) 
Theorem \ref{thm:main}, together with $\alpha_k = 1/\sqrt{k}$ and $\beta_k = \beta^k$, guarantees that, for all $K \geq 1$ and all $\bm{x}\in \mathbb{R}^d$,
\begin{align}\label{ineq:3_1}
\begin{split}
\frac{1}{K} \sum_{k=1}^K V_k(\bm{x})
&\leq
\frac{d D H}{2 \tilde{b}}\frac{s}{\alpha_K K}
+
\frac{G^2}{2 \tilde{b} \tilde{\gamma}^2 h_0^*}\frac{1}{s K} \sum_{k=1}^K \alpha_k
+
\frac{\sqrt{dD} G}{\tilde{b}} \frac{1}{K} \sum_{k=1}^K \beta^k\\
&\leq
\frac{d D H}{2 \tilde{b} \alpha}\frac{s}{\sqrt{K}}
+
\frac{G^2 \alpha}{\tilde{b} \tilde{\gamma}^2 h_0^*}\frac{1}{s \sqrt{K}}
+
\frac{\beta \sqrt{dD} G}{\tilde{b} \tilde{\beta}} \frac{1}{K},
\end{split}
\end{align}
where we use $\tilde{\beta} := 1 - \beta$, 
\begin{align*}
&\frac{1}{K} \sum_{k=1}^K \frac{1}{\sqrt{k}}
\leq 
\frac{1}{K} \left(1 + \int_1^K \frac{\mathrm{d} t}{\sqrt{t}}  \right)
= 
\frac{1}{K} \left(2 \sqrt{K} -1 \right)
\leq 
\frac{2}{\sqrt{K}},\\
&\frac{1}{K} \sum_{k=1}^K \beta^k
\leq  
\frac{1}{K} \sum_{k=1}^{+ \infty} \beta^k
= 
\frac{\beta}{\tilde{\beta}  K}.
\end{align*}
An argument similar to the one for showing (\ref{ineq:2}) and (\ref{ineq:3}) ensures that
(\ref{ineq:3_1}) implies that 
\begin{align}\label{ineq:3_3}
\min_{k\in [K]} \mathbb{E}\left[ \|\nabla f (\bm{x}_k)\|^2 \right]
\leq
\underbrace{\frac{d D H}{2 \tilde{b}\alpha}}_{A_\alpha} \frac{s}{\sqrt{K}}
+
\underbrace{\frac{G^2 \alpha}{\tilde{b} \tilde{\gamma}^2 h_0^*}}_{B_\alpha} \frac{1}{s \sqrt{K}} 
+
\underbrace{\frac{\beta \sqrt{dD} G}{\tilde{b} \tilde{\beta}}}_{C_\beta} \frac{1}{K}.
\end{align} 

(ii)
A sufficient condition for (\ref{evaluation}), i.e., 
\begin{align*}
\min_{k\in [K]} \mathbb{E}\left[ \|\nabla f (\bm{x}_k)\|^2 \right]
\leq \epsilon^2
\end{align*} 
is that the right-hand side of (\ref{ineq:3_3}) is equal to $\epsilon^2$, i.e., 
\begin{align*}
A_\alpha s^2 \sqrt{K} + B_\alpha \sqrt{K} + C_\beta s - \epsilon^2 s K  = 0,
\end{align*}
which implies that 
\begin{align*}
K(s) = \left\{ \frac{(A_\alpha s^2 + B_\alpha) + \sqrt{(A_\alpha s^2 + B_\alpha)^2 + 4 \epsilon^2 C_\beta s^2}}{2 \epsilon^2 s} \right\}^2.
\end{align*}
We have that
\begin{align*}
\frac{\mathrm{d} \sqrt{K(s)}}{\mathrm{d}s}
=
\frac{(A_\alpha s^2 + B_\alpha) + \sqrt{(A_\alpha s^2 + B_\alpha)^2 + 4 \epsilon^2 C_\beta s^2}}{2 \epsilon^2 s^2 \sqrt{(A_\alpha s^2 + B_\alpha)^2 + 4 \epsilon^2 C_\beta s^2}}
(A_\alpha s^2  - B_\alpha) 
\begin{cases}
\displaystyle{< 0 \quad\text{ if } s \in ( 0, s^\star ),}\\ 
\displaystyle{= 0 \quad\text{ if } s = s^\star = \sqrt{\frac{B_\alpha}{A_\alpha}},}\\
\displaystyle{> 0 \quad\text{ if } s \in (s^\star, + \infty),}
\end{cases}
\end{align*}
which implies that $K(s)$ attains the minimum $K(s^\star)$ when $s = s^\star$.
$\Box$

\end{document}